\documentclass[11pt,amsfonts]{amsart}

\newtheorem{theorem}{Theorem}[section]

\newtheorem{corollary}[theorem]{Corollary}

\newtheorem{proposition}[theorem]{Proposition}

\newtheorem{remark}[theorem]{Remark}

\newtheorem{example}[theorem]{Example}

\newtheorem{conjecture}[theorem]{Conjecture}

\def\endproof{\qed \medskip}
\def\blacksquare{\hbox to .60em{\vrule width .60em height .60em}}

\begin{document}

\title[]{Topics in Conformally Compact Einstein Metrics}

\author[]{Michael T. Anderson}

\thanks{Partially supported by NSF Grant DMS 0305865}

\maketitle

\setcounter{section}{0}

\section{Introduction.}
\setcounter{equation}{0}

 Conformal compactifications of Einstein metrics were introduced by 
Penrose [38], as a means to study the behavior of gravitational fields 
at infinity, i.e. the asymptotic behavior of solutions to the vacuum 
Einstein equations at null infinity. This has remained a very active area 
of research, cf. [27], [19] for recent surveys. In the context of Riemannian 
metrics, the modern study of conformally compact Einstein metrics began with the 
work of Fefferman-Graham [26], in connection with their study of 
conformal invariants of Riemannian metrics. Recent mathematical work in 
this area has been significantly influenced by the AdS/CFT (or 
gravity-gauge) correspondence in string theory, introduced by Maldacena 
[36]. We will only comment briefly here on aspects of the AdS/CFT 
correspondence, and refer to [2], [42], [7] for general surveys. 

 In this paper, we discuss recent mathematical progress in this area, 
focusing mainly on global aspects of conformally compact Einstein metrics 
and the global existence question for the Dirichlet problem. One reason for 
this is that it now appears that the beginnings of a general existence theory 
for such metrics may be emerging, at least in dimension 4. Of course to date there 
is no general theory for the existence of complete Einstein metrics on 
manifolds, with two notable exceptions; the existence theory for 
K\"ahler-Einstein metrics due to Calabi, Yau, Aubin and others, and 
the existence theory in dimension 3, due to Perelman, Hamilton and Thurston. 
In contrast to the situation for compact 4-manifolds, an existence theory 
for conformally compact Einstein metrics may not be that far beyond the current 
horizon. 

  We discuss numerous open problems on this topic; some new results are also 
presented, cf. in particular Theorem 3.4 and the discussion and results in 
Sections 4 and 5. 

\medskip

 In brief, the contents of the paper are as follows. The groundwork is laid 
in \S 2, where we discuss the moduli space of conformally compact Einstein metrics 
and the boundary map to the space of conformal infinities. The general situation 
is also illustrated by the discussion of a simple but important class of examples, 
the static AdS black hole metrics. Section 3 deals with the general asympototic 
behavior of the metrics near conformal infinity, and the control of the asymptotic 
behavior by the metric at infinity. It will be seen that at least in even 
dimensions, this issue is now quite well understood. Then in Section 4 we turn 
to the analysis of the behavior of the metrics on compact regions, away from 
infinity, mostly in dimension 4 where the possible degenerations can be described 
in terms of orbifold and cusp degenerations. In Section 5, we conclude with a 
discussion of the possibility of actually finding examples where orbifold or 
cusp degenerations occur. 

\medskip

  I would like to thank David Calderbank, Tom Farrell, Lowell Jones, 
Claude LeBrun, Rafe Mazzeo and Michael Singer for discussions related 
to various issues in the paper. Thanks also to Vestislav Apostolov and 
collegues for organizing an interesting workshop at the CRM, Montreal 
in July, 04.

\section{Conformally compact Einstein metrics.}
\setcounter{equation}{0}

 Let $M$ be the interior of a compact $(n+1)$-dimensional manifold 
$\bar M$ with boundary $\partial M$. A complete Riemannian metric 
$g$ on $M$ is $C^{m,\alpha}$ conformally compact if there is a defining 
function $\rho$ on $\bar M$ such that the conformally equivalent 
metric
\begin{equation} \label{e2.1}
\widetilde g = \rho^{2}g 
\end{equation}
extends to a $C^{m,\alpha}$ metric on the compactification $\bar M$. 
Here $\rho$ is a smooth, non-negative function on $\bar M$ with 
$\rho^{-1}(0) = \partial M$ and $d\rho \neq 0$ on $\partial M$. The 
induced metric $\gamma  = \widetilde g|_{\partial M}$ is the boundary 
metric associated to the compactification $\widetilde g$. Since there 
are many possible defining functions, there are many conformal 
compactifications of a given metric $g$, and so only the conformal class 
$[\gamma]$ of $\gamma$ on $\partial M$, called conformal infinity, is 
uniquely determined by $(M, g)$. Clearly any manifold $M$ carries many 
conformally compact metrics but we are mainly concerned here with Einstein 
metrics $g$, normalized so that
\begin{equation} \label{e2.2}
Ric_{g} = -ng. 
\end{equation}
A simple computation for conformal changes of metric shows that if $g$ 
is at least $C^{2}$ conformally compact, then the sectional curvature 
$K_{g}$ of $g$ satisfies
\begin{equation} \label{e2.3}
|K_{g}+1| = O(\rho^{2}). 
\end{equation}
Thus, the local geometry of $(M, g)$ approaches that of hyperbolic space, 
and conformally compact Einstein metrics are frequently called 
asymptotically hyperbolic (AH), or also Poincar\'e-Einstein. All these 
notions will be used here interchangeably. The natural ``threshold 
level'' for smoothness is $C^{2}$, since even if $g$ is $C^{m,\alpha}$ 
conformally compact, $m > 2$, (2.3) cannot be improved to $|K_{g}+1| = 
o(\rho^{2})$ in general. 

 Mathematically, an obviously basic issue in this area is the Dirichlet 
problem for conformally compact Einstein metrics: given the topological 
data $(M, \partial M)$, and a conformal class $[\gamma]$ on $\partial M$, 
does there exist a conformally compact Einstein metric $g$ on $M$, 
with conformal infinity $[\gamma]$? In one form or another, this 
question is the basic leitmotiv throughout this paper. As will be 
seen later, uniqueness of solutions with a given conformal infinity 
fails in general. 

\medskip

 To set the stage, we first examine the structure of the moduli space 
of Poincar\'e-Einstein metrics on a given $(n+1)$-manifold $M$. Let 
$E^{m,\alpha}$ be the space of Poincar\'e-Einstein metrics on $M$ which 
admit a $C^{2}$ conformal compactification $\bar g$ as in (2.1), 
with $C^{m,\alpha}$ boundary metric $\gamma$ on $\partial M$. Here $0 
< \alpha < 1$, $m \geq 2$, and we allow $m = \infty$ or $m = \omega$, 
the latter corresponding to real-analytic. The topology on $E^{m,\alpha}$ 
is given by a weighted H\"older norm, cf. (2.7) below; briefly, the topology 
is somewhat stronger than the $C^{2}$ topology on metrics on $\bar M$ under a 
conformal compactification $\widetilde g$ as in (2.1). Let ${\mathcal E}^{m,\alpha} = 
E^{m,\alpha}/{\rm Diff}_{1}^{m+1,\alpha}(\bar M)$, where 
${\rm Diff}_{1}^{m+1,\alpha}(\bar M)$ is the group of $C^{m+1,\alpha}$ 
diffeomorphisms of $\bar M$ inducing the identity on $\partial M$, acting 
on $E$ in the usual way by pullback. Next, let $Met^{m,\alpha}(\partial M)$ 
be the space of $C^{m,\alpha}$ metrics on $\partial M$ and ${\mathcal C}^{m,\alpha} 
= {\mathcal C}^{m,\alpha}(\partial M)$ the corresponding space of 
pointwise conformal classes. 

 The natural boundary map, 
\begin{equation} \label{e2.4}
\Pi : {\mathcal E}^{m,\alpha} \rightarrow  {\mathcal C}^{m,\alpha}, \ \ 
\Pi [g] = [\gamma], 
\end{equation}
takes a conformally compact Einstein metric $g$ on $M$ to its conformal 
infinity on $\partial M$. Thus, global existence for the Dirichlet 
problem is equivalent to the surjectivity of $\Pi$, while uniqueness 
is equivalent to the injectivity of $\Pi$. 

 The following result describes the general structure of ${\mathcal E}$ 
and the map $\Pi$, building on previous work of Graham-Lee [29] and 
Biquard [15].
\begin{theorem} \label{t 2.1.}
{\rm  (Manifold structure [5], [6])} Let $M$ be a compact, oriented $(n+1)$-manifold 
with boundary $\partial M$ with $n \geq 3$. If ${\mathcal E}^{m,\alpha}$ is non-empty, 
then ${\mathcal E}^{m,\alpha}$ is a smooth infinite dimensional manifold. Further, 
the boundary map 
$$\Pi : {\mathcal E}^{m,\alpha} \rightarrow  {\mathcal C}^{m,\alpha} $$
is a $C^{\infty}$ smooth Fredholm map of index 0.
\end{theorem}

  When $ m < \infty$, ${\mathcal E}^{m,\alpha}$ has the structure of a 
Banach manifold, while ${\mathcal E}^{\infty}$ has the structure of a 
Fr\'echet manifold. For $n = 3$, one expects that Theorem 2.1 also holds for 
$m \geq 2$, but this is an open problem. 

\medskip

   Theorem 2.1 shows that the moduli space ${\mathcal E}$ has a very 
satisfactory global structure. In particular if $M$ carries some 
Poincar\'e-Einstein metric, then it also carries a large set of them, 
mapping under $\Pi$ to at least a variety of finite codimension in ${\mathcal C}$. 
Recall that a metric $g \in {\mathcal E}$ is a regular point of $\Pi$ if 
$D_{g}\Pi$ is surjective. Since $\Pi$ is Fredholm of index $0$, $D_{g}\Pi$ 
is injective at regular points; hence, by the inverse function theorem, 
$\Pi$ is a local diffeomorphism in a neighborhood of each regular point. 

\begin{remark} \label{r 2.2.}
 {\rm  Note that Theorem 2.1 does not hold, as stated, when $n = 1$, 
i.e. in dimension 2. In this case, the space ${\mathcal E}$ 
as defined above is infinite dimensional, but it becomes finite 
dimensional when one divides out by the larger group of diffeomorphisms 
isotopic to the identity on $\bar M$. This space of conformally compact 
(geometrically finite) hyperbolic metrics on a surface $\Sigma$ is a smooth, 
finite dimensional manifold, but the conformal infinity is unique. 
The boundary $\partial \Sigma$ is a collection of circles and there is 
only one conformal structure on $S^{1}$ up to diffeomorphism. 
In particular, $\Pi$ is not of index 0. 

  When $n = 2$, Einstein metrics are again hyperbolic, and the space of 
such metrics, modulo diffeomorphisms isotopic to the identity, is 
parametrized by the Teichm\"uller space of conformal classes on Riemann 
surfaces forming $\partial M$. Thus, Theorem 2.1 does hold for $n = 2$. 
However, we point out that the map $\Phi$ in (2.7) below used in constructing 
${\mathcal E}$ is not Fredholm when $n = 1, 2$. Thus, the proof of Theorem 2.1 
does not extend to the case $n = 2$. }
\end{remark}

\medskip

 It is worthwhile to examine the local structure of the boundary map 
$\Pi$ near singular points in more detail. To do this, we need to 
discuss some background material, related to the proof of Theorem 2.1. 
Given a boundary metric $\gamma$, one may form the standard 
``hyperbolic cone'' metric on $\gamma$ by setting, in a neighborhood of 
$\partial M$,
$$g_{\gamma} = \rho^{-2}(d\rho^{2} +\gamma ), $$
and extending $g_{\gamma}$ to $M$ in a fixed but arbitrary way. Given a 
fixed background metric $g_{0}\in E^{m,\alpha}$ with boundary metric 
$\gamma_{0}$, for $\gamma $ near $\gamma_{0}$, let $g(\gamma) = g_{0} 
+ \eta (g_{\gamma}$ - $g_{\gamma_{0}})$, where $\eta$ is a cutoff 
function supported near $\partial M$. Thus $g(\gamma )$ is close to 
$g_{0}$ and consider metrics $g$ near $g_{0}$ of the form 
\begin{equation} \label{e2.5}
g = g(\gamma ) + h, 
\end{equation}
where $h$ is a symmetric bilinear form on $M$ which decays as 
$O(\rho^{2})$. 

 Essentially following [15], the Bianchi-gauged Einstein 
operator at $g_{0}$ is defined by 
\begin{equation} \label{e2.6}
\Phi (g) = Ric_{g} + ng + \delta_{g}^{*}\beta_{g(\gamma )}(g). 
\end{equation}
We view $\Phi$ as a map
\begin{equation} \label{e2.7}
\Phi : Met^{m,\alpha}(\partial M)\times {\mathbb S}_{2}^{m,\alpha}(M) 
\rightarrow  {\mathbb S}_{2}^{m-2,\alpha}(M), 
\end{equation}
$$\Phi (\gamma , h) = \Phi (g(\gamma ) + h), $$
where ${\mathbb S}_{2}^{m,\alpha}(M)$ is the space of symmetric bilinear 
forms $h$ on $M$, of the form $h = \rho^{2}\bar h$, with $\bar h$ 
bounded in $C^{m,\alpha}(M)$. It turns out that if $g_{0}\in 
E^{m,\alpha}$ then the variety $\Phi^{-1}(0)$ forms a local slice for 
the action of diffeomorphisms on $E^{m,\alpha}$ near $g_{0}$. 

 The derivative of $\Phi$ at $g_{0}$ with respect to the second factor 
is the linearized Einstein operator 
\begin{equation} \label{e2.8}
L(h) = D^{*}Dh - 2R(h), 
\end{equation}
$h\in{\mathbb S}_{2}^{m,\alpha}(M).$ By [29], this map is Fredholm, and so 
has finite dimensional kernel and cokernel. Let $K$ be the kernel of 
$L$ on ${\mathbb S}_{2}^{m,\alpha}(M)$; $K$ is also the kernel of $L$ on 
$L^{2}(M, g)$. To prove Theorem 2.1, it suffices to show that $\Phi$ 
is a submersion at any $g_{0}\in E^{m,\alpha}$, and for this one needs 
to show that the pairing
\begin{equation} \label{e2.9}
\int_{M}\langle D\Phi^{\widetilde g}(\dot \gamma_{0}), \kappa \rangle 
dV_{g_{0}}
\end{equation}
is non-degenerate, in the sense that for any $\kappa  \in  K_{g_{0}}$, 
there exists a variation $\dot \gamma_{0}$ of $\gamma_{0} = \Pi (g_{0})$ 
such that (2.9) is non-zero. This is actually not so easy in general, and 
we refer to [6] for details. 

\medskip

 The boundary map $\Pi$ is locally, near $g_{0}$, just the projection 
map on the first factor of $\Phi^{-1}(0)$ in (2.7). Thus, locally, a 
slice for ${\mathcal E}^{m,\alpha}$ through $g_{0}$ is written as a 
(possibly multi-valued and singular) graph over 
$Met^{m,\alpha}(\partial M)$. The kernel $K$ of $D\Pi$ at $g$ is the 
subspace at which the graph is vertical, and corresponds to the kernel 
$K$ of the operator $L$ in (2.8). 

  To understand the singularities of $\Pi$ in more detail, 
note that since $\Pi$ is Fredholm, it is locally proper, i.e. 
for any $g\in{\mathcal E}^{m,\alpha}$, there exists an open set ${\mathcal U}$ 
with $g\in{\mathcal U}$ such that $\Pi|_{{\mathcal U}}$ is a proper map onto 
its image ${\mathcal V} \subset {\mathcal C}$. This means that $\Pi$ has a 
local degree, ${\rm deg}_{g}\Pi \in {\mathbb Z}$, cf. [41], [13]; in fact if 
${\mathcal U}$ is chosen sufficiently small, then ${\rm deg}_{g}\Pi = -1, 0$ 
or $+1$. If ${\rm deg}_{g}\Pi \neq 0$, then $\Pi$ is locally surjective onto 
a neighborhood of $\gamma  = \Pi(g)$; this may or not be the case if 
${\rm deg}_{g}\Pi= 0$. Observe however that (of course) ${\rm deg}_{g}\Pi$ 
is not continuous in $g$. 

 The local degree can be calculated by examining the behavior of $\Pi$ 
on generic, finite dimensional slices. Thus, let $B$ be any 
$p$-dimensional local affine subspace (or submanifold) of 
$Met^{m,\alpha}(\partial M)$ with $\gamma  = \Pi(g) \in  B$ and 
consider the restriction of $\Phi$ to $B\times {\mathbb S}_{2}^{m,\alpha}(M)$, 
and correspondingly, the graph $E_{B}^{m,\alpha} = 
\Phi^{-1}(0)\cap\Pi^{-1}(B)$ of $E^{m,\alpha}$ over $B$. For a generic 
choice of $B$, $E_{B}^{m,\alpha}$ is a $p$-dimensional manifold, and thus 
one can examine the behavior of $\Pi|_{E_{B}^{m,\alpha}}$ in the 
context of the study of singularities of smooth mappings between 
equidimensional manifolds. By construction, cf. [13] for instance, one 
has for generic $B$, 
$${\rm deg}_{g}\Pi = {\rm deg}_{g}\Pi|_{E_{B}^{m,\alpha}}. $$
Consider for example the situation where $B$ is 1-dimensional. Then 
$E_{B}^{m,\alpha}$ is a local curve in $Met^{m,\alpha}(\partial M)\times 
{\mathbb S}_{2}^{m,\alpha}(M)$ graphed over the interval 
$B = (-\varepsilon ,\varepsilon)$, with 0 corresponding to $\gamma$. 
One sees that if ${\rm deg}_{g}\Pi|_{E_{B}^{m,\alpha}} = \pm 1,$ then 
$\Pi$ is locally surjective near $\gamma ,$ while if 
${\rm deg}_{g}\Pi|_{E_{B}^{m,\alpha}} = 0$, then locally 
$\Pi|_{E_{B}^{m,\alpha}}$ is a fold map, equivalent to 
$x \rightarrow  x^{2}$ on $(-\varepsilon , \varepsilon)$. In this case, at 
least in a small neighborhood ${\mathcal U}$ of $g$, $\Pi$ is not 
surjective onto a neighborhood of $\gamma$; there is a local ``wall'' in 
${\mathcal C}$, (the image of the fold locus), which $\Pi({\mathcal U})$ 
does not cross. 

 Some natural questions related to this discussion are the following: 
is the set of critical points of $\Pi$ a non-degenerate critical 
submanifold (in sense of Bott)? Is it possible that $\Pi$ maps a 
connected manifold or variety of dimension $\geq 1$ onto a point 
$\gamma\in{\mathcal C}$? 

\medskip

 At this point, it is useful to illustrate the discussion on the basis 
of some concrete examples.
\begin{example} \label{ex2.3}{\rm (Static AdS black hole metrics). 
Let $N^{n-1}$ be any closed $(n-1)$-dimensional manifold, which carries 
an Einstein metric $g_{N}$ satisfying
\begin{equation} \label{e2.10}
Ric_{g_{N}} = k(n-2)g_{N}, 
\end{equation}
where $k = +1, 0$ or $-1$. We assume $n \geq 3$. Consider the metric 
$g_{m}$ on ${\mathbb R}^{2}\times N$ defined by
\begin{equation} \label{e2.11}
g_{m} = V^{-1}dr^{2} + Vd\theta^{2} + r^{2}g_{N}, 
\end{equation}
where 
\begin{equation} \label{e2.12}
V(r) = k+r^{2}-\frac{2m}{r^{n-2}}. 
\end{equation}
Here $r \in  [r_{+},\infty)$, where $r_{+}$ is the largest root of $V$, 
and the circular parameter $\theta \in [0,\beta]$, where 
\begin{equation} \label{e2.13}
\beta  = 4\pi r_{+}/(nr_{+}^{2}+k(n-2)). 
\end{equation}
This choice of $\beta$ is required so that the metric $g_{m}$ is smooth 
at the locus $\{r = r_{+}\}$; if $\beta$ is arbitrary, the metric will 
have cone singularities normal to the locus $\{r = r_{+}\}$, although 
the metric is otherwise smooth. Since this locus is the fixed point set 
of the isometric $S^{1}$ action given by rotation in $\theta$, the set 
$\{r = r_{+}\}$ is diffeomorphic to $N$ and is totally geodesic; it 
corresponds to the horizon of the black hole.  A simple computation 
shows that the metrics $g_{m}$ are Einstein, satisfying (2.2). Further, 
it is easy to see these metrics are smoothly conformally compact; the 
conformal infinity of $g_{m}$ is given by the conformal class of the 
product metric on $S^{1}(\beta )\times (N, g_{N})$. 

 We discuss the cases $k > 0$, $k = 0$, $k < 0$ in turn.

{\bf I.}  Suppose $k = +1$. 

 As a function of $m \in  (0,\infty)$, observe that $\beta$ has a 
maximum value of $\beta_{0} = 2\pi (\frac{n-2}{n})^{1/2}$, and for 
every $m \neq  m_{0}$, there are two values $m^{\pm}$ of $m$ giving the 
same value of $\beta$. Thus two metrics have the same conformal 
infinity; in particular, the boundary map $\Pi$ in (2.4) is not 1-1 
along this curve. This behavior is the first example of non-uniqueness 
for the Dirichlet problem, and was discovered in [32] in the context of 
the AdS Schwarzschild metrics, where $N = S^{2}(1)$. 

 The map $\Pi$ is a fold map, (of the form $x \rightarrow  x^{2}$), in 
a neighborhood of the curve $g_{m}$ near $m = m_{0}$. The local degree 
at $g_{m_{0}}$ is 0 and $\Pi$ is not locally surjective. In fact, Theorem 2.4 
below implies that $\Pi$ is globally not surjective, in that the conformal 
class of $S^{1}(L)\times (N, g_{N})$, for $L > \beta_{0}$, is not in 
Im $\Pi$, cf. [5]. Observe that this result requires global smoothness 
of the Einstein metrics; if one allows cone singularities along the 
horizon $N = \{r = r_{+}\}$, i.e. if $\beta$ is allowed to be 
arbitrary, then one can go past the ``wall'' through 
$S^{1}(\beta_{0})\times (N, g_{N})$. This clearly illustrates the global 
nature of the global existence or surjectivity problem.

{\bf II.} Suppose $k = 0$.

 In this case $\beta  = 4\pi r_{+}/(nr_{+}^{2})$ is a monotone function 
of $r_{+}$ or $m$, so that it assumes all values in ${\mathbb R}^{+}$ as $m 
\in  (0,\infty)$. On the curve $g_{m}$, $\Pi$ is 1-1. 

 However, the actual situation is somewhat more subtle than this. 
Suppose for instance that $N = T^{n-1}$, so that $M = {\mathbb 
R}^{2}\times T^{n-1}$ is a solid torus. Topologically, the disc $D^{2} = 
{\mathbb R}^{2}$ can be attached onto {\it any} simple closed curve in the 
boundary $\partial M = T^{n}$ instead of just the ``trivial'' $S^{1}$ 
factor in the product $T^{n} = S^{1}\times T^{n-1}$. The resulting manifolds 
are all diffeomorphic. This can also be done metrically, preserving the 
Einstein condition, cf. [4], and leads to the existence of infinitely many 
distinct Einstein metrics on ${\mathbb R}^{2}\times T^{n-1}$ with the same 
conformal infinity $(T^{n}, [g_{0}])$, where $g_{0}$ is any flat metric. 

 Each of these metrics lies in a distinct component of the moduli space 
${\mathcal E}$, so that ${\mathcal E}$ has infinitely many components. This 
situation is closely related to the mapping class group $SL(n, {\mathbb Z})$ 
of $T^{n}$, i.e. the group of diffeomorphisms of $T^{n}$ modulo those 
homotopic to the identity map, (so called ``large diffeomorphisms''). 
Any element of $SL(n, {\mathbb Z})$ extends to a diffeomorphism of the 
solid torus ${\mathbb R}^{2}\times T^{n-1}$, and while $SL(n, {\mathbb Z})$ 
acts trivially on the moduli space of flat metrics on $T^{n}$, the action 
on ${\mathcal E}$ is highly non-trivial, giving rise to the distinct 
components of ${\mathcal E}$. Similar constructions can obviously 
be carried out for manifolds $N$ of the form $N = T^{k}\times N'$, $k \geq 1$, 
but it would be interesting to investigate the most general version of this 
phenomenon.

{\bf III.}  Suppose $k = -1$. 

 Again $\beta$ is a monotone function of $m$, and so takes on all values 
in ${\mathbb R}^{+};$ the boundary map $\Pi$ is 1-1 on the curve $g_{m}$. 
Further aspects of this case are discussed later in \S 5. }
\end{example}

 These simple examples already show a number of subtle features of the 
global behavior of the boundary map $\Pi$. With regard to the global 
surjectivity question, the basic property that one needs to make 
progress is to understand whether $\Pi$ is a proper map; if $\Pi$ is 
not proper, it is important to understand exactly what possible 
degenerations of Poincar\'e-Einstein metrics can or do  occur with 
controlled conformal infinity. Recall that $\Pi$ is proper if and only 
if $\Pi^{-1}(K)$ is compact in ${\mathcal E}$, whenever $K$ is compact in 
${\mathcal C}$. 

 If $\Pi$ is proper, then one has a well-defined ${\mathbb Z}_{2}$-valued 
degree, cf. [41]. In fact, since the spaces ${\mathcal E}$ and ${\mathcal C}$ 
can be given a well-defined orientation, one has a ${\mathbb Z}$-valued degree, 
given by
\begin{equation} \label{e2.14}
{\rm deg} \ \Pi = \sum_{g_{i}\in\Pi^{-1}[\gamma ]}(-1)^{ind_{g_{i}}}, 
\end{equation}
where $[\gamma]$ is a regular value of $\Pi$ and $ind_{g_{i}}$ is the 
$L^{2}$ index of $D_{g_{i}}\Pi$, i.e. the number of negative 
eigenvalues of the operator $L$ in (2.8) at $g_{i}$ acting on $L^{2}$, 
cf. [5]. Of course if deg $\Pi \neq 0$, then $\Pi$ is surjective; (if 
deg $\Pi = 0$, then $\Pi$ may or may not be surjective). Note that deg $\Pi$ is 
defined on each component ${\mathcal E}_{0}$ of ${\mathcal E} $ and may differ 
on different components. 

 Let $M = M^{4}$ be a 4-manifold, satisfying 
\begin{equation} \label{e2.15}
H_{2}(\partial M, {\mathbb R}) \rightarrow H_{2}(M, {\mathbb R}) \rightarrow 0.
\end{equation}
It is proved in [5] that $\Pi$ is then proper, when restricted to the space 
${\mathcal E}^{0}$ of Einstein metrics whose conformal infinity is of non-negative 
scalar curvature. More precisely, 
\begin{equation} \label{e2.16}
\Pi^{0}: {\mathcal E}^{0} \rightarrow  {\mathcal C}^{0} 
\end{equation}
is proper, where ${\mathcal C}^{0}$ is the space of conformal classes 
having a {\it non-flat} representative of non-negative scalar curvature and 
${\mathcal E}^{0} = \Pi^{-1}({\mathcal C}^{0})$; in particular there are only 
finitely many components to ${\mathcal E}_{0}$, compare with Example 2.3, 
Case II above. 

 In situations where $\Pi$ is proper, the degree can be calculated in 
a number of concrete situations by the following:

\begin{theorem} \label{t 2.4.} {\rm (Isometry Extension, [5])}
  Let $(M^{n+1},$ g) be a $C^{2}$ conformally compact Einstein metric 
with $C^{\infty}$ boundary metric $\gamma$, $n \geq 3$. Then any 
connected group $G$ of conformal isometries of $(\partial M, \gamma)$ 
extends to a group $G$ of isometries of $(M, g)$. 
\end{theorem}

 This result has a number of immediate consequences. For instance, it 
implies that the Poincar\'e (or hyperbolic) metric is the unique $C^{2}$ 
conformally compact Einstein metric on an $(n+1)$-manifold with conformal 
infinity given by the round metric on $S^{n}$; see also [12], [39] for previous 
special cases of this result. In particular, one has on $(B^{4}, S^{3})$, 
$${\rm deg} \ \Pi^{0} = 1, $$
so that $\Pi$ is surjective onto ${\mathcal C}^{0}$. On the other hand, on 
$(M^{4}, S^{3}), M^{4} \neq  B^{4}$, 
$${\rm deg} \ \Pi^{0} = 0, $$
since $\Pi$ cannot be surjective in this case. Another application of 
Theorem 2.4 is the following:

\begin{corollary} \label{c 2.5.}
  Let $M$ be any compact $(n+1)$-manifold with boundary $\partial M$, 
$n \geq 3$, and let $\hat M = M\cup_{\partial M}M$ be the closed manifold 
obtained by doubling $M$ across its boundary. Suppose $\partial M$ admits 
an effective $S^{1}$ action, but $\hat M$ admits no effective 
$S^{1}$ action. Then $\Pi = \Pi(M)$ is not surjective; in fact 
$${\rm Im} \Pi \cap  Met_{S^{1}}(\partial M) = \emptyset  , $$
where $Met_{S^{1}}(\partial M)$ is the space of $S^{1}$ invariant 
metrics on $\partial M.$ The space $Met_{S^{1}}(\partial M)$ is of 
infinite dimension and codimension in $Met(\partial M)$. 
\end{corollary}
{\endproof}

 As a simple example, let $\hat M = \Sigma_{g}\times N$, where $\Sigma$ 
is any surface of genus $g \geq 1$ and $N$ is any $K(\pi ,1)$ 
manifold with $\pi$ having no center; e.g. $N$ has a metric of 
non-positive curvature. Let $\sigma$ be a closed curve in $\Sigma_{g}$ 
which disconnects $\Sigma_{g}$ into two diffeomorphic components 
$\Sigma^{+}$ and $\Sigma^{-}$ with common boundary $\sigma$, and let 
$M = \Sigma^{+}\times N$. By [22], $\hat M$ does not admit an effective 
$S^{1}$ action, but of course $\partial M = S^{1}\times N$ admits such 
actions. Hence, Corollary 2.5 holds for such $M$. 

  On the other hand, if $\Sigma = S^{2}$ is of genus 0, then $M = 
{\mathbb R}^{2}\times N$ does admit $S^{1}$-invariant Poincar\'e-Einstein 
metrics, as discussed in Example 2.3. 

\medskip

 A basic issue is to extend the theory described above beyond boundary 
metrics of non-negative scalar curvature ${\mathcal C}^{0}$. This will be 
one of the themes discussed below. We begin with the analysis of 
Poincar\'e-Einstein metrics near the boundary, i.e. conformal infinity.

\section{Behavior near the Boundary.}
\setcounter{equation}{0}

 In this section, we study the behavior of Poincar\'e-Einstein metrics in 
a neighborhood of conformal infinity $(\partial M, \gamma)$. 

 For many purposes, the most natural compactifications are those 
defined by geodesic defining functions. Thus, a compactification 
$\bar g = \rho^{2}g$ as in (2.1) is called geodesic if $\rho (x) = 
dist_{\bar g}(x, \partial M)$. Each choice of boundary metric 
$\gamma \in [\gamma]$ determines a unique geodesic defining function 
$\rho$. For a geodesic compactification, one typically loses one 
derivative in the possible smoothness, but this will not be of major 
concern here, cf. also [11, App.B] on restoring loss of derivatives. 

 The Gauss Lemma gives the splitting
\begin{equation} \label{e3.1}
\bar g = d\rho^{2} + g_{\rho}, \ \ g = \rho^{-2}(d\rho^{2} + 
g_{\rho}), 
\end{equation}
where $g_{\rho}$ is a curve of metrics on $\partial M$. A simple and 
natural idea to examine the behavior of $g$ near infinity is to expand 
the curve of metrics $g_{\rho}$ on $\partial M$ in a Taylor series in 
$\rho .$ Surprisingly (at first), this turns out not always to be 
possible, as discovered in [26]. It turns out that the exact form of 
the expansion depends on whether $n$ is odd or even. If $n$ is odd, i.e 
$M$ is even-dimensional, then 
\begin{equation} \label{e3.2}
g_{\rho} \sim  g_{(0)} + \rho^{2}g_{(2)} + .... + \rho^{n-1}g_{(n-1)} + 
\rho^{n}g_{(n)} + \rho^{n+1}g_{(n+1)} + ... 
\end{equation}
This expansion is even in powers of $\rho$ up to order $n$. The 
coefficients $g_{(2k)}$, $2k \leq (n-1)$ are locally determined via the 
Einstein equations (2.2) by the boundary metric $\gamma  = g_{(0)}$. 
They are explicitly computable expressions in the curvature of $\gamma 
$ and its covariant derivatives, although their complexity grows 
rapidly with $k$. The term $g_{(n)}$ is transverse-traceless, i.e.
\begin{equation} \label{e3.3}
tr_{\gamma}g_{(n)} = 0,  \ \ \delta_{\gamma}g_{(n)} = 0, 
\end{equation}
but is otherwise undetermined by $\gamma$ and the Einstein equations; it 
depends on the particular structure of the AH Einstein metric $(M, g)$ near 
infinity. If $n$ is even, one has
\begin{equation} \label{e3.4}
g_{\rho} \sim  g_{(0)} + \rho^{2}g_{(2)} + .... + \rho^{n-2}g_{(n-2)} + 
\rho^{n}g_{(n)} + \rho^{n}\log\rho \ {\mathcal H}  + \rho^{n+1}g_{(n+1)} + 
... 
\end{equation}
Again (via the Einstein equations) the terms $g_{(2k)}$ up to order $n-2$ 
are explicitly computable from the boundary metric $\gamma$, as is the 
coefficient ${\mathcal H}$ of the first $\log\rho$ term. The term 
${\mathcal H}$ is transverse-traceless. The term $g_{(n)}$ satisfies 
\begin{equation} \label{e3.5}
tr_{\gamma}g_{(n)} = \tau , \ \  \delta_{\gamma}g_{(n)} = \delta ,  
\end{equation}
where again $\tau$ and $\delta$ are explicitly determined by the 
boundary metric $\gamma$ and its derivatives; however, as before 
$g_{(n)}$ is otherwise undetermined by $\gamma$. There are $(\log\rho)^{k}$ 
terms that appear in the expansion at order $> n$. 

  Note also that these expansions (3.2) and (3.4) depend on the choice 
of boundary metric. Transformation properties of the coefficients $g_{(i)}$, 
$i \leq n$, under conformal changes have been explicitly studied in the physics 
literature, cf. [24]. As discovered by Fefferman-Graham [26], the term 
${\mathcal H}$ is conformally invariant, or more precisely covariant: 
if $\widetilde \gamma = \phi^{2}\gamma$, then $\widetilde {\mathcal H} = 
\phi^{2-n} {\mathcal H}$. 

\begin{remark} \label{r3.1}
{\rm Analogous to the Fefferman-Graham expansion above, there is a 
formal expansion of a vacuum solution to the Einstein equations near 
null infinity, although this has been carried out in detail only in 
dimension 3+1, cf. [16]. This expansion is closely related to the 
properties of the Penrose conformal compactification. More recently, 
as discussed in [20], logarithmic terms appear in the expansion 
in general, and these play an important role in understanding the 
global structure of the space-time. }

\end{remark}

 Mathematically, it is of some importance to keep in mind that the 
expansions (3.2), (3.4) are only formal, obtained by conformally 
compactifiying the Einstein equations and taking iterated Lie derivatives 
of $\bar g$ at $\rho = 0$;
\begin{equation} \label{e3.6}
g_{(k)} = \frac{1}{k!}{\mathcal L}_{T}^{(k)}\bar g, 
\end{equation}
where $T = \nabla\rho$. If $\bar g \in  C^{m,\alpha}(\bar M)$, 
then the expansions hold up to order $m+\alpha$. However, boundary 
regularity results are needed to ensure that if an AH Einstein metric 
$g$ with boundary metric $\gamma$ satisfies $\gamma  \in  
C^{m,\alpha}(\partial M)$, then the compactification $\bar g \in  
C^{m,\alpha}(\bar M)$ or $C^{m',\alpha'}(\bar M)$. 

 In both cases $n$ odd or even, the Einstein equations determine all 
higher order coefficients $g_{(k)}$ (and coefficients of the $\log$ 
terms), in terms of $g_{(0)}$ and $g_{(n)}$, so that an AH Einstein 
metric is formally determined by $g_{(0)}$ and $g_{(n)}$ near $\partial M$. 
The term $g_{(0)}$ corresponds to Dirichlet boundary data on $\partial M$, 
while $g_{(n)}$ corresponds to Neumann boundary data, (in analogy with the 
scalar Laplace operator). Thus, on AH Einstein metrics, the formal 
correspondence
\begin{equation} \label{e3.7}
g_{(0)} \rightarrow  g_{(n)} 
\end{equation}
is analogous to the Dirichlet-to-Neumann map for harmonic functions. 
However, the map (3.7) is only well-defined if there is a unique AH 
Einstein metric with boundary data $\gamma  = g_{(0)}$; as seen above 
on the curve of AdS Schwarzschild metrics for example, this is not 
always the case. Understanding the correspondence (3.7) is a basic 
issue, both mathematically and in certain aspects of the AdS/CFT 
correspondence. Again in a formal sense, knowing $g_{(0)}$ and 
$g_{(n)}$ allows one to locally construct the bulk gravitational field, 
i.e. the Poincar\'e-Einstein metric, at least near $\partial M$ via the 
expansion (3.2) or (3.4). 

\medskip

 To begin to make some of the discussion above more rigorous, we next 
discuss the boundary regularity issue; many aspects of this have been 
resolved over the past few years. Suppose first $n = 3$, so dim $M  = 
4$. If $g \in  E^{m,\alpha}$, $m \geq 2$, then by definition $g$ has a 
$C^{2}$ conformal compactification to a $C^{m,\alpha}$ boundary metric 
$\gamma$. In [4], it is proved that there is a $C^{m,\alpha}$ 
conformal compactification $\widetilde g \in  C^{m,\alpha}(\bar M)$ 
of $g$, cf. also [6]. This result also holds if $m = \infty$ or $m = 
\omega$. It is proved using the fact that 4-dimensional Einstein 
metrics satisfy the Bach equations, cf. [14], which are conformally 
invariant. In suitable gauges, the Bach equation can be recast as a 
non-degenerate elliptic system of equations for a conformal compactification 
$\widetilde g$, and the result follows from elliptic boundary regularity. 

 In dimension 4, the Bach tensor is the Fefferman-Graham obstruction 
tensor ${\mathcal H}$ above. In any even dimension, the system of 
equations 
\begin{equation} \label{e3.8}
{\mathcal H}  = 0 
\end{equation}
is conformally invariant, and is satisfied by metrics conformal to 
Einstein metrics. Thus, one might expect that the method using the Bach 
equation in [4], [6] when $n = 3$ can be extended to all $n$ odd. 
This is in fact the case, and has been worked out in detail by 
Helliwell [31]. Thus, essentially the same regularity results hold 
for $n$ odd. 

 When $n$ is even, so that dim $M$ is odd, this type of boundary 
regularity cannot hold of course, due to the presence of the 
logarithmic terms in the FG expansion. A result of Lee [35] shows that 
if $g \in  E^{m,\alpha}$ and $m < n$, then $g$ is $C^{m,\alpha}$ 
conformally compact. This is optimal, but does not reach the important 
threshold level $m = n$, where logarithmic terms and the important 
$g_{(n)}$ term first appear. Recently, Chru\'sciel et al. [21] have 
proved that when $g \in E^{\infty}$, i.e. $g$ has a $C^{\infty}$ 
boundary metric $\gamma$, then $g$ has a $C^{\infty}$ 
polyhomogeneous conformal compactification, so that the expansion (3.4) 
exists as an asymptotic series. Moreover, if $\gamma \in 
C^{m,\alpha}(\partial M)$, then the expansion exists up to order $k$, 
where $k$ can be made large by choosing $m$ sufficiently large; 
(in general $m$ must be much larger than $k$). Finally, it has recently 
been proved by Kichenassamy [34] that when $g \in E^{\omega}$ and 
$g_{(n)}$ is real-analytic, the formal series (3.4) exists, i.e. it 
is summable, and it converges to $g_{\rho}$.

\medskip

 These results have the following immediate consequence. Suppose $n$ is 
odd. Given any real-analytic symmetric bilinear forms $g_{(0)}$ and 
$g_{(n)}$ on $\partial M$, satisfying (3.3), there exists a unique 
$C^{\omega}$ conformally compact Einstein metric $g$ defined in a 
thickening $\partial M\times [0,\varepsilon)$ of $\partial M$. If 
instead $n$ is even, given any analytic symmetric bilinear forms $g_{(0)}$ 
and $g_{(n)}$ on $\partial M,$ satisfying (3.5), there exists a unique 
$C^{\infty}$ polyhomogeneous conformally compact Einstein metric $g$ 
defined in a thickening $\partial M\times [0,\varepsilon )$ of $\partial M$. 
In both cases, the expansions (3.2) or (3.4) converge to the metric 
$g_{\rho}$. These results follow from the work in [4], [6], [31] when $n$ 
is odd, and [34] when $n$ is even. Since analytic data $g_{(0)}$ and 
$g_{(n)}$ may be specified arbitrarily and independently of each other, 
subject only to the constraint (3.3) or (3.5), to give ``local'' AH 
Einstein metrics, defined in a neighborhood of $\partial M$, this shows 
that the correspondence (3.7) must depend highly on global properties 
of Poincar\'e-Einstein metrics.

 On the other hand, it is well-known that the use of analytic data to 
solve elliptic-type problems is misleading. While the Dirichlet or 
Neumann problem is formally well-posed, the Cauchy problem is not. 
Standard examples involving Laplace operator and harmonic functions 
show that even if Cauchy data on a boundary converge smoothly to limit 
Cauchy data on the boundary, the corresponding solutions do not converge 
to a limit in any neighborhood of the boundary. 

 To pass from analytic to smooth boundary data, one needs apriori 
estimates or equivalently a stability result. In this respect, one has 
the following:

\begin{theorem} \label{t 3.2} {\rm (Local Stability, [4], [6])}
  Let $g$ be a $C^{2}$ conformally compact Einstein metric, defined in 
a region $\Omega  = [0,\rho_{0}]\times \partial M$ containing $\partial M$, 
where $\rho$ is a geodesic compactification. Suppose there exists a 
compactification $\widetilde g$ with $C^{m,\alpha}$ boundary metric 
$\gamma$, such that
\begin{equation} \label{e3.9}
||\widetilde g||_{C^{1,\alpha}(\Omega)} \leq  K. 
\end{equation}
If $n = 3$, then there is a (possibly different) compactification, 
also called $\widetilde g$, such that, in $\Omega' = [0,\frac{\rho_{0}}{2}]
\times \partial M$, one has the estimate 
\begin{equation} \label{e3.10}
||\widetilde g||_{C^{m,\alpha}(\Omega')} \leq  C, 
\end{equation}
where the constant $C$ depends only on $K$, $m$, $\alpha$, $n$ and 
$\rho_{0}$. 
\end{theorem}

  This result is proved simultaneously with the boundary regularity result 
itself, i.e. using the fact that $\widetilde g$ is a solution of the Bach 
equations together with standard estimates for solutions of elliptic systems 
of PDE's. 

  Using similar ideas as discussed above in connection with (3.8), Theorem 3.2 
also holds for all $n$ odd, at least if $C^{2}$ is replaced by $C^{n,\alpha}$ 
and $C^{1,\alpha}$ is replaced by $C^{n,\alpha}$ in (3.9), with $m > n$ in 
(3.10), cf. [31]. 

\medskip

 Theorem 3.2 shows that if two solutions are close in a weak norm, 
($C^{1,\alpha}$ or $C^{n,\alpha}$), then they are close in a strong norm, 
$C^{m,\alpha}$, $m > n$, provided the boundary metrics are 
close in a strong norm. It would be very interesting if a similar 
result can be proved when $n$ is even, (i.e. in odd dimensions). A 
direct generalization is of course not possible, due again to the 
logarithmic terms. Redefining the H\"older norms to take such logarithmic 
terms into account, it would be very surprising if such a stability 
result did not hold; however, a proof remains to be established. 

 In even dimensions, the local stability theorem allows one to pass to 
limits in the analytic data problem above. Thus, suppose $\gamma  = 
g_{(0)}$ and $g_{(n)}$ are arbitrary $C^{m,\alpha}$ data on $\partial M$, 
subject to the constraint (3.3). Let $\gamma_{i}$ and 
$(g_{(n)})_{i}$ be a sequence of analytic data satisfying (3.3) 
converging to $\gamma$ and $g_{(n)}$ in the $C^{m,\alpha}$ topology, 
(such sequences always exist), and let $\widetilde g_{i}$ be the 
corresponding sequence of conformal compactifications of 
Poincar\'e-Einstein metrics defined in regions $\Omega_{i}$. If the 
hypothesis (3.9), (with $C^{1,\alpha}$ replaced by $C^{n,\alpha}$ for 
$n > 3$), held on the sequence $\{\widetilde g_{i}\}$, i.e. 
$\Omega_{i} = \Omega$ is uniform, then it follows that 
$\{\widetilde g_{i}\}$ converges in the $C^{m,\alpha}$ topology on 
$\Omega$ to a limit $\widetilde g \in  C^{m,\alpha}(\Omega )$. The metric 
$\widetilde g$ is a conformal compactification of a Poincar\'e-Einstein 
metric $g$, defined at least on $\Omega$. In other words, it would then 
follow that arbitrary smooth $\gamma$ and $g_{(n)}$ can be realized as 
local boundary data.

 However, the following result shows this cannot be the case:

\begin{theorem} \label{t 3.3.} {\rm (Unique Continuation, [8])}
   Let data $(g_{(0)}, g_{(n)})$ be arbitrarily given, satisfying the 
constraints (3.3) or (3.5), in some open set $U \subset \partial M$, 
with $(g_{(0)}, g_{(n)}) \in C^{m,\alpha}(U)$, for $m > n$ and any 
$n \geq 3$. If $g$ is a $C^{m,\alpha}$ conformally compact Einstein 
metric realizing the data $(g_{(0)}, g_{(n)})$, defined in a neighborhood 
$\Omega$ with $\Omega\cap\partial M = U$, then $g$ is the unique such 
metric, up to local isometry. 
\end{theorem}

 This result implies in particular that local Cauchy data in an open 
set $U \subset \partial M$ determine the global behavior of the 
metric, and the topology of the manifold, up to covering spaces; here 
we use the fact that Einstein metrics are real-analytic in the 
interior, and so trivially satisfy a unique continuation property. 
It follows that $(g_{(0)}, g_{(n)})$ in $U$ necessarily determine 
$(g_{(0)}, g_{(n)})$ outside $U$. (This is of course obvious for analytic 
data $(g_{(0)}, g_{(n)})$ on $\partial M$).

 It then follows that, for $\{\widetilde g_{i}\}$ as above, the weak 
uniform bound (3.9) cannot hold in general. The metrics must degenerate 
in a small neighborhood $\Omega$ of $U\subset \partial M$, for ``most'' 
choices of $g_{(n)}$, given any fixed choice of $\gamma = g_{(0)}$ on 
$\partial M$. 

\medskip

 We now contrast this situation with the situation for {\it globally} 
defined Poincar\'e-Einstein metrics. For emphasis, for the result below 
we require that $(M, g)$ is globally conformally compact, i.e. $M$ is the 
interior of a compact manifold with boundary, and $g$ is complete and 
globally defined on $M$. 

\begin{theorem} \label{t 3.4.} {\rm  (Control near Boundary)}
 Let $(M^{n+1}, g)$ be a globally conformally compact Poincar\'e-Einstein 
metric, with $n$ odd, so that dim $M$ is even. Suppose that $g$ is 
$C^{2}$ conformally compact, with $C^{m,\alpha}$ boundary metric 
$\gamma$, with $m > n$ and $m \geq 6$ if $n = 3$. Then there exists a 
neighborhood $\Omega = [0,\rho_{0}]\times \partial M$ of $\partial M$, 
depending only on the boundary data $(\partial M, \gamma)$ such that
\begin{equation} \label{e3.11}
||\widetilde g||_{C^{m,\alpha}(\Omega)} \leq  C, 
\end{equation}
in some compactification $\widetilde g$.
\end{theorem}

 The bound (3.11) implies that the boundary map $\Pi$ is proper near 
conformal infinity, in the sense that if one has a fixed boundary 
metric $\gamma$, or compact set of boundary metrics $\gamma \in  
\Gamma$, then the set of Poincar\'e-Einstein metrics with boundary 
metric $\gamma$, (or $\gamma\in\Gamma$), is compact, as far as their 
behavior in $\Omega$ is concerned; any sequence has a convergent 
subsequence on a fixed domain $\Omega$, where $\Omega$ only depends 
on the boundary data. 

\noindent
{\bf Proof:}
 This result is proved for $n = 3$ in [5], and the proof for arbitrary 
$n$ odd is very similar. Thus, we refer to [5] for much of the proof, 
and only discuss those situations where the proof needs to be modified 
in higher dimensions. 

 There are several steps in the proof. First, let $\bar g$ be the 
geodesic compactification of $g$ determined by $\gamma$, and let $\tau$ 
be the distance to the cutlocus of the normal exponential map from 
$(\partial M, \gamma)$ into $(M, \bar g)$. Here of course $g$ is 
any Poincar\'e-Einstein metric on $M$ with boundary metric $\gamma$, (or 
$\gamma\in\Gamma$). The first (and most important) step is to prove that 
there is a constant $\tau_{0} > 0$, depending only on $n$ and $\gamma$ (or 
$\Gamma$) such that
\begin{equation} \label{e3.12}
\tau (x) \geq  \tau_{0}. 
\end{equation}
The estimate (3.12) already implies for instance that the topology of 
$M$ cannot become non-trivial too close to the boundary $\partial M$. 
The proof of (3.12) in [5, Prop.4.5] holds with only minor and essentially 
obvious changes in all even dimensions, given the local stability result, 
Theorem 3.2. As noted in [5, Remark 2.4], one should use the renormalized action 
in place of the renormalized volume or $L^{2}$ norm of the Weyl curvature. Also, 
the classification of ${\mathbb R}^{n}$-invariant solutions as AdS toral 
black holes is given in [9], (again the proof of this holds in all dimensions). 

 Next, let $\zeta (x) = \zeta^{n,\alpha}(x)$ be the $C^{n,\alpha}$ harmonic 
radius of $(M, \bar g)$ at $x$, for a fixed $\alpha < 1$. The next claim, 
(cf. [5, Prop.4.4]) is that there is a constant $\zeta_{0}$, depending only 
on $n$ and $\gamma$, such that
\begin{equation} \label{e3.13}
\zeta (x) \geq  \zeta_{0}\tau (x). 
\end{equation}
(The proof in [5, Prop.4.4] uses the $L^{p}$ curvature radius, but the 
proof works equally well for the much stronger $C^{n,\alpha}$ harmonic 
radius). 

  The proof of (3.13) is by contradiction. If (3.13) does not hold, then 
there exist $x_{i}\in (M_{i}, g_{i})$ such that $\xi (x_{i}) << \tau 
(x_{i})$. Choose $x_{i}$ to realize the minimum of the ratio $\xi /\tau$. 
One then takes a blow-up limit of the rescalings $g_{i}'  = \zeta 
(x_{i})^{-2}\bar g_{i}$ based at $x_{i}$. Since $\zeta'(x_{i}) = 1$, 
$\zeta'(y_{i}) \geq \frac{1}{2}$ for $y_{i}$ within bounded 
$g_{i}'$-distance to $x_{i}$. It follows that in a subsequence, one has 
convergence to a complete limit $(N, g', x)$. The local stability result, 
Theorem 3.2, implies that the convergence to the limit is in the (strong) 
$C^{n,\alpha}$ topology. The radius $\zeta$ is continuous in this topology, 
and hence the limit $(N, g')$ cannot be flat, since $\zeta'(x) = 1$. Here 
one must also use the non-collapse or volume comparison estimates in 
[5, Lemma3.8ff]. Thus, to obtain a contradiction, it suffices to prove 
that the limit $(N, g')$ must be flat. To do this, one distinguishes the 
following two situations:

 I. $dist_{g_{i}'}(x_{i}, \partial M_{i}) \leq D$, for some $D <  
\infty$. In this case, the limit $N$ has a boundary $(\partial N, 
\gamma')$. Since this limit is the blow-up of $(\partial M_{i}, \gamma_{i})$, 
it is clear that $(\partial N, \gamma')$ is flat $({\mathbb R}^{n}, 
\delta)$, where $\delta$ is the flat metric. (Here we use of course 
the fact that $\Gamma$ is compact). Moreover, $\partial N$ is totally 
geodesic in $N$. As in [5], $(N, g')$ is Ricci-flat, $Ric_{g'} = 0$. 
The proof that $N$ is actually flat in [5] used the fact that $N$ contains 
a line; when $n = 3$, i.e. in dimension 4, this implies $N$ is flat. This 
of course does not hold in higher dimensions. Instead, one can argue as 
follows. Since $N$ is Ricci-flat and has flat and totally geodesic boundary 
${\mathbb R}^{n}$, the reflection double of $N$ across ${\mathbb R}^{n}$ is a 
weak $C^{1}$ solution of the Einstein equations $Ric_{g'} = 0$ on 
${\mathbb R}^{n+1}$. Elliptic regularity implies that $g'$ is then real-analytic 
across ${\mathbb R}^{n}$. It then follows easily from the Cauchy-Kovalevsky 
theorem that $(N, g')$ is flat. This gives the required contradiction in this case. 

 II. $dist_{g_{i}'}(x_{i}, \partial M_{i}) \rightarrow \infty$ as $i 
\rightarrow  \infty$. For this case, we give a different and simpler 
proof than that in [5, Prop.4.4, Case II]. 

 Let $d_{i}(x) = dist_{g_{i}'}(x_{i}, \partial M_{i})$. It follows 
easily from Case I above that
$$\zeta_{i}(y_{i}) \geq  (1-\delta )d_{i}(y_{i}),  $$
for $y_{i}$ within bounded distance to $(\partial M, g_{i}')$, with 
$\delta \rightarrow 1$ as $i \rightarrow \infty$. This just 
corresponds to the statement that the geometry becomes flat near 
$\partial M_{i}$ with respect to $g_{i}'$, which has been proved in 
Case I. Now by hypothesis, at $x_{i}$, $\zeta_{i}(x_{i})/d_{i}(x_{i}) 
\rightarrow 0$, (since the ratio is scale-invariant and 
$\zeta_{i}(x_{i}) = 1$ in the scale $g_{i}'$). Therefore, by continuity, 
there are points $y_{i}$ such that $\zeta (y_{i}) = \frac{1}{2}d(y_{i})$, 
with $\zeta (z_{i}) \geq  \frac{1}{2}d(z_{i})$, for all $z_{i}$ such that 
$d(z_{i}) \leq  d(y_{i})$. One now works in the scale $\hat g_{i} = 
\zeta_{i}(y_{i})^{-2}g_{i}$ where $\hat \zeta_{i}(y_{i}) = 1$ and hence 
$dist_{\hat g_{i}}(y_{i}, \partial M) = 2$. The proof is now completed just 
as in Case I. Thus, one may pass to a limit $(N, \hat g, y)$. On the one 
hand, the limit $(N, \hat g)$ is not flat, since, by Theorem 3.2 the 
convergence to the limit is in $C^{n,\alpha}$ and $\zeta$ is continuous 
in this topology, so that $\hat \zeta(y) = 1$. As before, $(N, \hat g)$ has 
flat and totally geodesic boundary, and the same proof as in Case I implies 
that $(N, \hat g)$ is flat, giving again a contradiction. 

 Taken together, (3.12) and (3.13) imply that $\zeta(x) \geq \tau_{1} > 0$, 
for all $x$ in a neighborhood of $\partial M$ of fixed size in $(M, \bar g)$. 
The bound (3.11) is then a consequence of the local stability result, 
Theorem 3.2. 
{\endproof}

 An odd dimensional analogue of Theorem 3.4 is unknown, and it would be 
very interesting to know if a suitable version of it holds. The exact 
formulation would of course have to be modified somewhat, due to the 
logarithmic terms. Using Kichenassamy's result [34], Javaheri [33] has 
proved an analogue of Theorem 3.4 in odd dimensions in the context of 
analytic boundary metrics.

\section{Behavior away from the Boundary.}
\setcounter{equation}{0}

 At least in even dimensions, the analysis in \S 3 shows that the 
global behavior of the boundary map $\Pi$ depends only on the behavior 
of Einstein metrics in the interior, a fixed distance away from the 
boundary, (depending only on the boundary metric), in a geodesic 
compactification. Thus, the issue of whether $\Pi$ is proper becomes a 
question on the behavior of Einstein metrics in the interior, i.e. on 
compact sets, away from infinity; the structure near infinity is 
uniformly controlled by the data at infinity. 

 Thus, in effect, one is dealing with the behavior of Einstein metrics 
on compact manifolds (with boundary). Presumably, the degeneration of 
such metrics has the same general features as the degeneration of 
Einstein metrics on compact manifolds without boundary. A detailed 
study of degenerations of Einstein metrics on compact 4-manifolds was 
first carried out in [3]. Since there is no general theory of such 
degenerations in higher dimensions, we restrict in this section to 
dimension 4. 

 Let $\{g_{i}\}$ be a sequence of Poincar\'e-Einstein metrics on a fixed 
4-manifold $M$, with conformal infinities $\{\gamma_{i}\}\subset\Gamma$, 
where $\Gamma$ is a compact set in ${\mathcal C}^{m,\alpha}$. There are 
three possibilities for the behavior of $\{g_{i}\}$, in subsequences; 
cf. [5] for a more detailed discussion.

 I. {\em Convergence}: A subsequence of $\{g_{i}\}$ converges, modulo 
diffeomorphisms, to a limit Poincar\'e-Einstein metric $g$ on $M$, with 
boundary metric $\gamma\in\Gamma$. There is a compactification 
$\widetilde g_{i} = \rho^{2}g_{i}$ of $g_{i}$ such that the subsequence 
$\{\widetilde g_{i}\}$ converges in the $C^{m,\alpha}$ topology on $\bar M$. 

 II. {\em Orbifolds}: A subsequence of $\{g_{i}\}$ converges, modulo 
diffeomorphisms, to a limit Poincar\'e-Einstein orbifold-singular metric 
$g$ on $M$, with boundary metric $\gamma\in\Gamma$. The singular metric 
$g$ is a smooth metric on an orbifold $V$, and $M$ is a smooth resolution 
of $V$. There are only a finite number of singularities, each the vertex 
of a cone on a spherical space form. Away from the singularities, the 
convergence is smooth, as in I. The subsequence $(M, g_{i})$ converges to 
$(V, g)$ in the Gromov-Hausdorff topology, [30]. 

 III. {\em Cusps}: A subsequence of $\{g_{i}\}$ converges, modulo 
diffeomorphisms, to a limit Poincar\'e-Einstein metric with cusps $g$ 
on a connected manifold $N$, with boundary metric $\gamma\in\Gamma$, 
possibly with a finite number of orbifold singularities.

 More precisely, the limit $(N, g)$ has conformal infinity $(\partial M, 
\gamma)$, but has in addition a collection of complete, finite volume 
ends. Thus there is a compact hypersurface $H \subset N$, 
disconnecting $N$ into two non-compact regions, the outside and inside; 
the outside region contains conformal infinity, and so has the same number 
of components as $\partial M$, while the inside region is connected and the 
metric $g$ is complete and of finite volume. 

 If $x_{i}$ are base points in $(M, g_{i})$ within bounded distance to 
$\partial M$ in the geodesic compactification $\bar g_{i}$, then 
the sequence $(M, g_{i}, x_{i})$ converges in the {\it pointed} 
Gromov-Hausdorff topology, (cf. [30]), to the limit $(N, g, x)$. The 
convergence is smooth, in the sense of I, away from any orbifold singular 
points, and uniform on compact sets in $\bar N = N \cup \partial N$, where 
$\partial N = \partial M$. Note that if one chooses other base points $y_{i}$ 
in $(M, g_{i})$ with $dist_{\bar g_{i}}(y_{i}, \partial M) \rightarrow \infty$, 
then $(M, g_{i}, y_{i})$ may limit on other complete, finite volume 
manifolds $(N', g', y)$; see the discussion regarding the case of 
surfaces below. However, since they play no role in the analysis here, 
we ignore these other components of the limit. 

\medskip

 If the boundary map $\Pi$ is to be proper, one must show that only 
the convergence case above occurs, i.e. rule out the possible formation 
of orbifold and cusp limits. We discuss these in turn. 

 The orbifold limits are topological, in the sense that essential 
2-cycles in $M$ not coming from $\partial M$ must be collapsed to 
points under $\{g_{i}\}$. Thus, for example, if one has a surjection
\begin{equation} \label{e4.1}
H_{2}(\partial M, {\mathbb R} ) \rightarrow  H_{2}(M, {\mathbb R} ) 
\rightarrow  0, 
\end{equation}
for instance $M = B^{4}$, then orbifold limits cannot occur, cf. [5]. The 
condition (4.1) however is not necessary, and there are 4-manifolds not 
satisfying (4.1) which do not admit any orbifold degenerations. A list of 
some examples is given in [5, \S 6]. Also, it appears that the families of 
self-dual Poincar\'e-Einstein metrics constructed by Calderbank-Singer [17], 
which are natural analogues of the Gibbons-Hawking metrics, do not admit 
orbifold degenerations, [18]. In fact, at this time, it seems that the only 
known example of orbifold degenerations is that of the Taub-Bolt curve of 
metrics on $T(S^{2})$ which degenerate to the orbifold $C({\mathbb R}{\mathbb P}^{3})$, 
cf. [5]; I am grateful to Michael Singer for pointing this out. 

\begin{remark} \label{r4.1}
{\rm In this context, it is worth noting that the manifold theorem, 
Theorem 2.1, holds also for orbifolds. Thus, let $V$ be an $(n+1)$-dimensional 
orbifold with boundary, in the sense that $V$ is a smooth manifold away from
from finitely many singular points in the interior, each having a neighborhood 
homeomorphic to the cone on a spherical space form. (Note that this definition 
of orbifold is much more restrictive than the general definition due to 
Thurston). Let ${\mathcal E}^{m,\alpha}(V)$ be the moduli space of orbifold smooth 
Poincar\'e-Einstein metrics on $V$, defined as in \S 2. Then Theorem 2.1 holds 
for ${\mathcal E}^{m,\alpha}(V)$; the proof is exactly the same. In fact, all 
the discussion and results above, from \S 2 to the classification I-III above, 
holds equally well for Poincar\'e-Einstein orbifold metrics. 

  As will be seen in the following, it would be very interesting to understand 
to what extent Theorem 2.1 generalizes to metrics with other singularities (e.g. 
cusps) on a compact manifold with boundary, cf. also [37] for some further 
results in this direction. }
\end{remark}

 Although the situation of orbifold degenerations still needs to be better 
understood in general, the issue of cusp formation is much more serious and much 
less well-understood. As seen in Example 2.3(II), there are at least some 
situations where cusp degenerations can occur. There are no known relations 
between the possibility of cusp formation and the topology of $M$, (as is 
the case with orbifold degenerations); this is a fundamental and very interesting 
open problem, which exists also for Einstein metrics on compact manifolds. 
It would also be useful to obtain more detailed information about 
the geometry of cusp ends.

\medskip

 Instead of trying to find situations where orbifold and cusp formation 
can be ruled out, (as in (4.1) for example), one can take a different 
perspective. Namely, these are the only possible degenerations of 
Poincar\'e-Einstein metrics with controlled conformal infinity and so 
it is natural to consider an enlarged space of Poincar\'e-Einstein 
metrics which includes these limits. 

 Thus, let $\bar{\mathcal E}$ be the completion of the moduli space 
${\mathcal E}$ of Poincar\'e-Einstein metrics with respect to the pointed 
Gromov-Hausdorff topology; the base points $x$ are 
chosen so that 
\begin{equation} \label{e4.2}
dist_{\bar g}(x, \partial M) = 1, 
\end{equation}
for example. The discussion concerning I-III above shows that metrics in 
$\bar{\mathcal E}$ with controlled conformal infinity $\gamma \in \Gamma$ are 
compact in this topology; any sequence in $\bar{\mathcal E}$ has a convergent 
subsequence to a limit in $\bar{\mathcal E}$ with conformal infinity 
$[\gamma] \in \Gamma$. Note that this topology on $\bar{\mathcal E}$ is quite 
different than the (unpointed) Gromov-Hausdorff topology; if $\{g_{i}\}$ is a 
sequence in ${\mathcal E}$ converging to a cusp metric in $\bar{\mathcal E}$, 
then $dist_{GH}(g_{i}, g_{0}) \rightarrow  \infty$, for any fixed 
$g_{0}\in{\mathcal E}$. This is because $diam_{g_{i}}M \rightarrow \infty$, 
and the diameter is continuous in the Gromov-Hausdorff topology. In 
particular, although the Gromov-Hausdorff topology is a metric 
topology on ${\mathcal E}$, this is not known for the pointed Gromov-Hausdorff 
topology on $\bar{\mathcal E}$. 

 Let $\partial{\mathcal E} = \bar{\mathcal E}\setminus {\mathcal E}$, 
so that $\partial{\mathcal E}$ consists of orbifold and cusped 
Poincar\'e-Einstein metrics, obtained as limits of smooth 
Poincar\'e-Einstein metrics on $M$. If $g_{i}\in{\mathcal E}$ converges to 
$g\in\partial{\mathcal E}$, then for any fixed $R$, the metrics $g_{i}$ on 
$B_{x_{i}}(R)$ converge smoothly to the limit metric $g$ on $B_{x}(R)$, 
away from any orbifold singular points; here $x_{i}$ are base points 
satisfying (4.2) and $x_{i} \rightarrow x$. Briefly, away from 
orbifold singular points, one has smooth convergence on compact subsets. 
Further, the compactified metrics $\bar g_{i}$ converge in $C^{m,\alpha}$ 
to $\bar g$ up the boundary $\partial M$. (As always, the smooth convergence 
is understood to be modulo diffeomorphisms). Note that the closure of 
${\mathcal E}$ in the Gromov-Hausdorff topology consists of ${\mathcal E}$ 
together with orbifold-singular Poincar\'e-Einstein metrics obtained as 
limits.

 Now one has an extension $\bar{\Pi}$ of $\Pi$ to $\bar{\mathcal E}$, 
and 
\begin{equation} \label{e4.3}
\bar{\Pi}: \bar{\mathcal E} \rightarrow  {\mathcal C}  
\end{equation}
is continuous, cf. [5]. Moreover, by construction, $\bar{\Pi}$ is 
proper. 

 If $\bar{\mathcal E}$ has roughly the structure of a manifold, then 
as is the case with ${\mathcal E}^{0}$ before, one can define a degree 
${\rm deg} \ \bar{\Pi}$ associated with each component of $\bar{\mathcal E}$
and 
$${\rm deg} \ \bar{\Pi} = {\rm deg} \ \Pi. $$
So ${\rm deg} \  \bar{\Pi} \neq 0$ implies at least that almost every choice 
of conformal class in ${\mathcal C}$ is the conformal infinity of a smooth 
Poincar\'e-Einstein metric on $M$.  

 Unfortunately, very little is known about the structure of 
$\bar{\mathcal E}$, even regarding its point set topology. As a first 
step, the following conjecture seems very plausible:
\begin{conjecture} \label{c4.2}
{\rm  For any component ${\mathcal E}_{0}$ of ${\mathcal E}$, 
$\bar \Pi(\partial {\mathcal E}_{0})$ has empty interior in 
${\mathcal C}$.  }
\end{conjecture}

  The intuition leading to Conjecture 4.2 is that $\Pi$ has Fredholm index 
0. As a simple illustration, let $E = \{(x,y,z)\in {\mathbb R}^{3}: z \in 
(0,1)\}$. Then $\bar E = E \cup \partial E$, $\partial E = \{z = 1\}$ is 
a manifold with boundary and the projection map $\pi: E \rightarrow 
{\mathbb R}^{2}$, $\pi(x,y,z) = (x,y,0)$ has the property that $\pi$ 
maps $\partial E$ onto $\pi(E)$. Of course Conjecture 4.2 fails on this 
example; however, the index of the map $\pi$ is one. 

  Similarly, if $\Pi(\partial {\mathcal E}_{0})$ did have non-trivial interior 
${\mathcal V} \subset {\mathcal C}$, and if $\overline{\mathcal E}_{0}$ is 
reasonably well-behaved, one would expect there are curves $\sigma$ in 
$\overline{\mathcal E}_{0}$ on which $\Pi$ is constant, i.e. for all 
$\gamma \in {\mathcal U}$, there exists $\sigma_{\gamma}(t) \subset 
\overline{\mathcal E}_{0}$, with $\sigma_{\gamma}(t) \subset 
{\mathcal E}_{0}$ for $t > 0$ and $\sigma_{\gamma}(0) \in 
\partial{\mathcal E}_{0}$, such that $\Pi \circ \sigma_{\gamma} = \gamma$. 
Hence, if $\gamma$ is a regular value of $\Pi$, then ${\rm index} \ \Pi \geq 1$, 
which is impossible. 

  As will be seen in \S 5, (cf. Example 5.2), Conjecture 4.2 is probably 
false if ${\mathcal E}_{0}$ is replaced by ${\mathcal E}$. As a toy model 
where this conjecture fails, (with ${\rm index} \ \pi = 0$), let $E_{i}$ 
be the collection of planes in ${\mathbb R}^{3}$ given by $E_{i} = 
\{z = 1 - \frac{1}{i}\}$. Now connect these planes by a collection of tubes 
or wormholes, deleting the corresponding discs in $\{E_{i}\}$ and let $E$ be 
the resulting connected space. Then $\bar E = E \cup E_{\infty}$, where 
$E_{\infty} = \{z = 1\}$. As above, let $\pi: E \rightarrow {\mathbb R}^{2}$ 
be the projection onto $E_{0} = \{z = 0\}$. One may then choose the connecting 
tubes so that $\pi$ is continuous and surjective on $E$, on $\bar E$, and on 
$\partial E = E_{\infty}$. By choosing the tubes to become arbitrarily small 
and dense near $\partial E$, one may arrange that $E$ is uniformly locally 
path connected. 

  One would not expect that ${\mathcal E}_{0}$ or $\overline{\mathcal E}_{0}$ 
has such a complicated structure. Instead, it seems more likely that both 
$\partial{\mathcal E}_{0}$ and $\Pi(\partial{\mathcal E}_{0})$ should be 
lower-dimensional in the spaces $\bar{\mathcal E}$ and ${\mathcal C}$ respectively. 
If codim $\partial{\mathcal E}_{0} = 1$ in $\bar{\mathcal E}_{0}$, then 
$\partial{\mathcal E}_{0}$ acts as a topological boundary and so 
$\bar{\mathcal E}_{0}$ does not have the structure of a manifold; at best 
it is a manifold with boundary. In this case, it will be difficult to define 
a suitable degree. On the other hand, if codim $\partial{\mathcal E}_{0} > 1$, 
then the metric boundary $\partial{\mathcal E}_{0}$ is not topological and one 
expects that $\bar{\mathcal E}_{0}$ behaves sufficiently well to allow one to 
define a degree deg $\bar{\Pi}$ on $\bar{\mathcal E}_{0}$. It would of course 
be very interesting to make progress on these speculative remarks. 

\medskip

 An interesting alternate path is to try extend the map $\Phi$ in 
(2.6)-(2.7) to singular metrics which effectively model orbifold 
singular and cusp metrics on the manifold $M$. This is perhaps 
easier in the orbifold case, since the behavior of the Einstein metrics 
$g_{i}$ converging to an orbifold limit $g$ is quite well-understood. 
If $\Phi$ can be extended to such an enlarged space, consisting of 
smooth and singular metrics on $M$ modelling orbifolds and cusps, such 
that $\Phi$ is still a smooth mapping, with Fredholm linearization $L$, 
then the same proof as Theorem 2.1 will show that $\bar{\mathcal E}$ is 
a smooth manifold. For a discussion of orbifold singular metrics on a 
manifold $M$, (as opposed to smooth metrics on an orbifold associated to 
$M$), cf. [3]. 

\medskip

 With regard to the work to follow in \S 5, it is worthwhile to 
describe in some detail the simplest situation where curves of Einstein 
metrics form cusps, i.e. the case of hyperbolic metrics on surfaces. 
Thus, let $(\Sigma, g_{\Sigma})$ be any complete conformally compact 
Riemann surface with non-empty boundary of constant negative curvature, 
normalized so that 
\begin{equation} \label{e4.4}
K_{\Sigma} = Ric_{\Sigma} = -1, 
\end{equation}
and with $\pi_{1}(\Sigma) \neq 0$. Topologically, $\Sigma$ is 
$S^{2}$ with at least two discs removed,  or a surface of genus $g \geq 1$, 
with at least one disc removed. The boundary $\partial\Sigma $ is 
a union of $q$ circles, $\partial\Sigma  = \cup_{1}^{q}S^{1}$, with 
$q \geq 1$. In the free homotopy class of each end $E_{i}$ of $\Sigma$, 
one has a unique closed geodesic $\sigma_{i}$, $1 \leq  i \leq q$, 
of length $\alpha_{i} > 0$. 

 Let ${\mathcal M}$ be the moduli space of such metrics satisfying (4.4). 
There are several definitions of the moduli space, depending on the 
choice of the action of the diffeomorphism group on $\partial\Sigma$. 
To obtain a finite dimensional space, ${\mathcal M}$ is considered as the 
space of all conformally compact metrics satisfying (4.4) divided out 
by the action of all diffeomorphisms of $\bar \Sigma = 
\Sigma\cup\partial\Sigma$ mapping $\partial\Sigma$ onto itself. It is 
well-known, cf. [1] for example, that ${\mathcal M}$ is a smooth 
orbifold, of dimension
\begin{equation} \label{e4.5}
m = dim_{{\mathbb R}}{\mathcal M} = 6g-6+3q, 
\end{equation}
(with $m = 2$ if $g = 0$ and $q = 2$). The boundary map
\begin{equation} \label{e4.6}
\Pi: {\mathcal M}  \rightarrow  {\mathcal C} /{\rm Diff} 
\end{equation}
is a constant map, since $S^{1}$ has a unique conformal structure up to 
diffeomorphism. Thus one has ${\rm index} \ \Pi = m > 0$. 

 The boundary $\partial{\mathcal M}$ of the moduli space ${\mathcal M}$ with 
respect to the Deligne-Mumford compactification consists of Riemann 
surfaces with nodes or punctures; this coincides with the boundary in 
the pointed Gromov-Hausdorff topology, where $\partial{\mathcal M}$ is 
represented by complete hyperbolic metrics which have cusp ends, 
obtained by shrinking a collection of disjoint closed geodesics in 
$\Sigma$ to 0 length. Note that such geodesics may or may not include 
geodesics from the collection $\{\sigma_{i}\}$. Thus $\partial{\mathcal M}$ 
is stratified by the moduli spaces of Riemann surfaces of lower genus, and 
a positive number of punctures; the strata are of dimension 
$$d = 6g-6+3q-2p, $$
where $p$ is the number of cusp ends, (punctures). In particular, 
$\partial{\mathcal M} $ has codimension 2 in $\overline{\mathcal M}$. The 
closure $\overline{\mathcal E} = {\mathcal M}\cup\partial{\mathcal M}$ 
has the structure of a real-analytic variety, cf. [1]. The boundary map 
extends to $\overline{\mathcal E}$, and it is still the constant map.

\medskip

 One would like to have a similar concrete description of cusp formation 
on some class of examples in higher dimensions. However, no such examples 
are known. Observe that conformally compact metrics are not closed under 
products; also the product of a compact metric and conformally compact 
metric is not conformally compact. In \S 5, we discuss a construction of 
families of conformally compact metrics forming cusps, based on the model 
of this behavior for surfaces. However, perhaps surprisingly, this does not 
lead to examples of Einstein metrics. 

 The following remains a simple but basic open question: does there 
exist a curve of Poincar\'e-Einstein metrics on $M^{n+1}$, $n \geq 3$, 
which converges to a Poincar\'e-Einstein metric with cusps? 

\begin{remark} \label{r4.4}
{\rm The classification of degenerations in I-III above is special to 
dimension 4 and very little is known in such generality in higher 
dimensions. However, in the presence of symmetry, the equations for 
Einstein metrics on higher dimensional manifolds can be reduced to 
the Einstein equations coupled to other fields in lower dimensions, via 
the well-known Kaluza-Klein procedure. In this regard, note that 
Theorem 2.4  implies that symmetries of a boundary metric $\gamma$ 
are automatically inherited by any Einstein metric $(N, g)$ filling 
$(\partial N, \gamma)$. 

  For example, suppose the compact group $G$ acts freely and isometrically 
on a Poincar\'e-Einstein metric $(N^{n+1}, g_{N})$. Let $M = N / G$ be the 
orbit space of this action; then the metric $g_{N}$ may be written in the 
form 
\begin{equation} \label{e4.7}
g_{N} = \pi^{*}g_{M} + \langle \theta, \theta \rangle,
\end{equation}
where $\pi: N \rightarrow M$ is the projection onto the orbit space, 
$\theta$ is a connection 1-form on $N$ with values in the Lie algebra 
${\mathcal L}(G)$ and $\langle, \rangle$ is a family of left-invariant 
metrics on ${\mathcal L}(G)$ parametrized by $x \in M$. 

   The Einstein equations (2.2) for $g_{N}$ become the Einstein equations 
for $g_{M}$ coupled to the gauge field $\theta$ and form $\langle, \rangle$. 
When dim $M = 4$, one can then consider whether the results above for the 
Einstein equations generalize to the Einstein equations coupled to 
various extra fields. This has been worked out in detail by Javaheri 
[33] for the case that $G = S^{1}$ and the action is static, so that 
the metric (4.7) has the form of a warped product; the equations on $M^{4}$ 
then take the form of the Einstein equations coupled to a scalar field. Note 
that already in this case, the Fefferman-Graham expansion on $M^{4}$ has 
logarithmic terms, due to the extra scalar field. }
\end{remark}

\section{Discussion on Cusp and Orbifold Formation.}
\setcounter{equation}{0}

 In this section, we show that it is not very easy to find continuous 
curves of Poincar\'e-Einstein metrics on a fixed manifold which limit on 
Poincar\'e-Einstein metrics with cusps. This gives some evidence, not 
particularly strong at the moment, but nevertheless suggesting that cusps 
may not form in components ${\mathcal E}_{0}$ of ${\mathcal E}$. Although 
the main focus of this section is on cusp formation, it will be seen that 
similar results often apply to orbifold formation. 

\medskip

 Let $(N, g_{0})$ be a Poincar\'e-Einstein metric with cusps. As in \S 2, 
let ${\mathbb S}_{2}^{m,\alpha}$ be the space of symmetric bilinear forms 
on $N$ which are bounded in $C^{m,\alpha}$ with respect to $g_{0}$ and 
decay in $C^{m,\alpha}$ at conformal infinity on the order of 
$\rho^{2}$. The map $\Phi  = \Phi^{g_{0}}$ is then defined as in 
(2.6)-(2.7). 
\begin{conjecture} \label{c 5.1.}
  Let (N, $g_{0})$ be a Poincar\'e-Einstein with cusps, and suppose that 
\begin{equation} \label{e5.1}
K_{N} \leq  0, 
\end{equation}
i.e. $g_{0}$ has non-positive curvature. Then the map $\Phi^{g_{0}}$ is 
a submersion at $g_{0}$, and the boundary map $\Pi$ taking 
Poincar\'e-Einstein metrics with cusps on $N$ to conformal classes ${\mathcal C}$ 
is a diffeomorphism in a neighborhood of $g_{0}$. 
\end{conjecture}

 Conjecture 5.1 stands in stark contrast to Conjecture 4.2; they 
almost contradict each other. In fact, they would contradict each other 
if one knew that Conjecture 4.2 holds and that every $\gamma\in {\rm Im} \ \Pi$ 
near $\gamma_{0}$ from Proposition 5.1 is the limit of a sequence 
$\gamma_{i} = \Pi(g_{i})$ with $g_{i}$ in a connected component 
${\mathcal E}_{0}$. If this were the case, it would follow that cusps 
satisfying (5.1) cannot form as limits within ${\mathcal E}_{0}$. Exactly 
the same remarks apply to orbifolds in place of cusps, since Conjecture 5.1 
does hold for orbifold Poincar\'e-Einstein metrics on an orbifold $V$.

\medskip

 In this context, it is worth considering some concrete examples:

\begin{example} \label{ex5.2}
{\rm  Let $g_{C}$ be the standard hyperbolic cusp metric on $N = {\mathbb R} 
\times T^{n}$ given by 
\begin{equation} \label{e5.2}
g_{C} = dr^{2} + e^{2r}\gamma_{T^{n}}, 
\end{equation}
where $\gamma_{T^{n}}$ is any flat metric on the torus $T^{n}$. Clearly 
$(N, g_{C})$ satisfies (5.1), so Conjecture 5.1 would imply that $\Pi$ is a 
local diffeomorphism near $g_{C}$, i.e. any boundary metric $\gamma$ 
near a flat metric $\gamma_{T^{n}}$ on $T^{n}$ is the conformal 
infinity of a complete Poincar\'e-Einstein cusp metric on ${\mathbb R} 
\times T^{n}.$ 

 On the other hand, as discussed in Example 2.3(II), there is an 
infinite sequence of conformally compact twisted toral black hole 
metrics $g_{i}$ on $M = {\mathbb R}^{2}\times T^{n-1}$. The metrics $g_{i}$ 
converge to $g_{C}$ in the pointed Gromov-Hausdorff topology, cf. [4]. These 
metrics all lie in distinct components ${\mathcal E}_{i}$ of ${\mathcal E}(M)$. 
If $\Pi_{i}: {\mathcal E}_{i} \rightarrow  {\mathcal C}$ is the boundary map, 
then we conjecture that for all $i$ large, $\Pi_{i}$ is surjective 
onto a fixed neighborhood ${\mathcal V}$ of $\Pi(g_{i}) = (T^{n}, \gamma_{T^{n}})$. 

 In this case, every conformal class $[\gamma] \in  {\mathcal V}$, for some open 
set ${\mathcal V} \subset {\mathcal C}$ containing $g_{T^{n}}$ is the 
conformal infinity of an infinite sequence of Poincar\'e-Einstein metrics on 
$M$, limiting on a Poincar\'e-Einstein cusp metric on $N$. This indicates that 
Conjecture 4.2 is false if the assumption that ${\mathcal E}_{0}$ is connected 
is dropped. 

 We point out that exactly the same discussion holds with ${\mathbb R} 
\times T^{n}$ replaced by any conformally compact hyperbolic manifold $N$, 
with a collection of cusp ends. As shown in [23], the cusp ends can be Dehn 
filled with solid tori to produce Poincar\'e-Einstein metrics $(M_{i}, 
g_{i})$ with a fixed conformal infinity. In this case, instead of 
having infinitely many components ${\mathcal E}_{i} = {\mathcal E}_{i}(M)$ 
of ${\mathcal E}$ on a fixed manifold $M$, one has a collection of components 
${\mathcal E}_{i} = {\mathcal E} (M_{i})$ on infinitely many topologically 
distinct manifolds $M_{i}$, with common boundary $\partial N$. }
\end{example}

 The discussion above presents some speculative evidence that cusps 
do not form within $\bar{\mathcal E}_{0}$, for any component 
${\mathcal E}_{0}$ of ${\mathcal E}$. Next, we present a construction of 
(connected) families of conformally compact metrics which are very close 
to being Einstein, and which do limit on cusps. This seems to be the simplest 
possible construction of such metrics, since it is based on the 
formation of cusps on surfaces. However, we argue that rather 
surprisingly, it is unlikely that these metrics can be perturbed to 
nearby Poincar\'e-Einstein metrics. 

\medskip

 To begin the construction, return to the static AdS black hole metrics 
(2.11) on ${\mathbb R}^{2}\times N,$ with $k = -1$. In this situation, 
$g_{m}$ is well-defined for negative values of $m$; in fact, $g_{m}$ is 
well-defined for $m \in  [m_{-}, \infty ),$ where 
\begin{equation} \label{e5.3}
m_{-} = -\frac{1}{n-2}(\frac{n-2}{n})^{n/2}, \ \ {\rm with} \ \  r_{+} = 
(\frac{n-2}{n})^{1/2}. 
\end{equation}
For the extremal value $m_{-}$ of $m$, $V(r_{+}) = V' (r_{+}) = 0$, and a 
simple calculation, (cf. (5.6) below) shows that the horizon $\{r = r_{+}\}$ 
occurs at infinite distance to any given point in $({\mathbb R}^{2}\times N, 
g_{m_{-}})$; the horizon in this case is called degenerate, (with zero surface 
gravity). Note that $\beta (m_{-}) = \infty$, so that the $\theta$-circles are 
in fact lines ${\mathbb R}$. As $m$ decreases to $m_{-}$, the horizon diverges 
to infinity, (in the opposite direction from the conformal infinity), while the 
length of the $\theta$-circles expands to $\infty$. Thus, the metric 
$g_{m_{-}}$ is a complete metric on the manifold ${\mathbb R}\times {\mathbb R}
\times N = {\mathbb R}^{2}\times N$, but is no longer conformally compact; 
the conformal infinity is ${\mathbb R}\times (N, g_{N})$. 

 However, one may divide the infinite $\theta$-factor ${\mathbb R}$ of the 
metric $g_{m_{-}}$ by ${\mathbb Z}$ to obtain a complete metric $g_{E}$ 
on $C\times N = {\mathbb R} \times S^{1}\times N$ of the form
\begin{equation} \label{e5.4}
g_{E} = V^{-1}dr^{2} + Vd\theta^{2} + r^{2}g_{N}, 
\end{equation}
where $V(r) = -1+r^{2}-\frac{2m_{-}}{r^{n-2}}$, with $m_{-}$ and 
$r_{+}$ given by (5.3). The length $\beta$ of the $\theta$-parameter 
in (5.4) is now arbitrary. The metric $g_{E}$ is called an extreme black 
hole metric, and is Poincar\'e-Einstein  with a single cusp end; $g_{E}$ 
has a smooth conformal compactification to the boundary metric 
$S^{1}(\beta)\times (N, g_{N})$. 

 To understand the behavior of the metric $g_{E}$ in the cusp region, 
we convert to geodesic coordinates.  Let $ds = V^{-1/2}dr$, so that, 
(up to an additive constant),
\begin{equation} \label{e5.5}
s = \int_{r_{+}+1}V^{-1/2}(r)dr.
\end{equation}
Thus, as $r \rightarrow  r_{+}$, $s \rightarrow  -\infty$, while as $r 
\rightarrow  \infty$, $s \rightarrow  \infty$. The metric (5.4) takes 
the form
\begin{equation} \label{e5.6}
g_{E} = ds^{2} + V(s)d\theta^{2} + r^{2}(s)g_{N}, 
\end{equation}
and the integral curves of $\nabla s$ are geodesics.  A simple 
calculation shows that
\begin{equation} \label{e5.7}
V(s) = ne^{2\sqrt{n}s}(1 + \varepsilon (s)), \ \ {\rm as} \ \  
s \rightarrow  -\infty , 
\end{equation}
where $\varepsilon (s) \rightarrow 0$ as $s \rightarrow -\infty$. 

 As $s \rightarrow  -\infty$, the length of the $\theta$-circles of 
course goes to 0, i.e. one has a collapse. However, the collapse can be 
unwrapped by passing to large covering spaces of the $S^{1}$ factor; on 
any sequence of base points $x_{i}$ with $s(x_{i}) \rightarrow  -\infty$, 
one may choose coverings so that the length of the $S^{1}$ factor at 
$x_{i}$ is approximately 1. One may then pass to a pointed limit to 
obtain the metric
\begin{equation} \label{e5.8}
g_{\infty} = ds^{2} + e^{2\sqrt{n}s}d\theta^{2} + r_{+}^{2}g_{N}. 
\end{equation}
The metric (5.8) is a product of the constant curvature metric on the 
cusp $C = {\mathbb R}\times S^{1}$ with a rescaling of the metric $g_{N}$ on 
$N$. By (2.10), the Ricci curvature of $r_{+}^{2}g_{N}$ equals 
$-r_{+}^{-2}(n-2) = -ng_{N}$, while the curvature of the cusp metric is 
also $-n$. The metric (5.8) (of course) has Ricci curvature 
$Ric_{g_{\infty}} = -ng_{\infty}$. The limit (5.8) is unique, up to 
rescalings of the length of the $S^{1}$ factor.

 The curvature of the extremal metric $g_{E}$ converges to that of 
$g_{\infty}$ exponentially fast in $s$. Straightforward computation from 
the estimates above shows that
\begin{equation} \label{e5.9}
||R_{g_{E}} - R_{g_{\infty}}|| \leq  Ce^{\sqrt{n}s}, \ \ {\rm as} \ \  
s \rightarrow  -\infty , 
\end{equation}
where the norm is the $L^{\infty}$ norm. The same estimate holds for 
all covariant derivatives of these curvatures.  

\medskip

 We now  make $g_{E}$ conformally compact, by closing off the 
cusp end. To do this, glue $\Sigma\times N$, where $\Sigma$ is any 
hyperbolic surface with an open expanding end, onto the cusp end of 
$g_{E}$. For simplicity, assume that $\Sigma$ has a single end; it is 
easy to generalize the construction below to any finite number of ends. 

 Thus, let $(\Sigma, g_{\Sigma})$ be any conformally compact Riemann 
surface with connected, non-empty boundary of constant negative 
curvature, normalized so that 
\begin{equation} \label{e5.10}
K_{\Sigma} = Ric_{\Sigma} = -n, 
\end{equation}
and with $\pi_{1}(\Sigma) \neq 0$. Let $\sigma$ be the unique closed 
geodesic in the free homotopy class of the end $E$ of $\Sigma$, and 
let $\alpha$ be the length of $\sigma$. We will only consider metrics 
$g_{\Sigma}$ in the moduli space ${\mathcal M}$, discussed in \S 4, for 
which
\begin{equation} \label{e5.11}
\alpha  \leq  \varepsilon , 
\end{equation}
where $\varepsilon$ is fixed and sufficiently small. Let ${\mathcal 
M}_{\varepsilon}$ be the domain in ${\mathcal M}$ satisfying (5.11). 

 Let $T(\tau)$ be the tubular neighborhood of radius $\tau$ about 
$\sigma$; coordinates may be introduced in this region so that the 
metric on $T(\tau)$ has the form
\begin{equation} \label{e5.12}
g_{\Sigma} = dt^{2} + \cosh^{2}(\sqrt{n}t)d\theta^{2}, 
\end{equation}
where the length of $\theta$ is $\alpha$ and $t \in  (-\tau , \tau)$. 
For $\alpha$ (arbitrarily) small, the expression (5.12) is valid 
for $\tau$ (arbitrarily) large.

 Next, on the product $\Sigma\times N$, form the product metric 
\begin{equation} \label{e5.13}
g_{D} = g_{\Sigma} + r_{+}^{2}g_{N}. 
\end{equation}
The metric $g_{D}$ is Einstein, of Ricci curvature $-n$.

 We now truncate the two metrics $g_{D}$ and $g_{E}$ and glue them 
together. To begin, topologically, set 
$$M = \Sigma\times N. $$
While the product metric $g_{D}$ on $M$ is Einstein, it is not 
conformally compact. This metric has one end $E$ of the form $C\times N$, 
where $C$ is an expanding cusp. In the $\tau$-tubular neighbhorhood of the 
geodesic $\sigma \subset E \subset \Sigma$, the metric has the form (5.12). 
Choose $R$ large (to be determined below), and let $D_{R}$ be the region in 
$\Sigma\times N$ where $t \leq R$, in the end $E$. Thus, $\partial D_{R} = 
S^{1}(L^{-})\times (r_{+}^{2}g_{N})$ where 
\begin{equation} \label{e5.14}
L^{-} = \cosh(\sqrt{n}R)\alpha . 
\end{equation}

 Next, take the conformally compact extreme metric $g_{E}$ on $C\times N$, 
and truncate it to the region $E^{R}$ where $s \geq -R$. The length of the 
boundary circle $\partial C$ is then 
\begin{equation} \label{e5.15}
L^{+} = V^{1/2}(-R)\beta . 
\end{equation}
To perform the glueing, we require that the lengths of the circles 
agree, $L^{-} = L^{+}$. Since the length $\beta$ is fixed, given the 
length $\alpha$, this imposes the relation
\begin{equation} \label{e5.16}
\alpha  = \{\cosh(\sqrt{n}R)\}^{-1}V^{1/2}(-R)\beta . 
\end{equation}
The parameter $t$ in (5.12) is related with the parameter $s$ in (5.6) 
by setting $t-R = s+R$. 

 Given these choices, one may easily construct a conformally compact 
approximate Einstein metric $\widetilde g$ on $M$ by attaching the 
truncated metrics $D_{R}$ and $E^{R}$ along their boundaries and 
smoothing the seam in a neighborhood $U$ of radius 1 of the boundaries. 
The metric $(M, \widetilde g)$ is smoothly conformally compact, and 
Einstein outside the glueing region $U$. Conformal infinity is given by 
the conformal class $(\partial M, [\gamma_{0}(\beta )])$, where 
$\partial M = S^{1}\times N$, and $\gamma_{0}(\beta) = \beta^{2}d\theta^{2} 
+ g_{N}$, with $\theta\in [0,1]$. 

  Let $\Phi  = \Phi^{\widetilde g}$ be as in (2.6)-(2.7). Then by construction, 
$\Phi^{\widetilde g}(\widetilde g) = 0$ outside $U$. An elementary computation, 
using (5.9) and simple estimates for the $2^{\rm nd}$ fundamental forms of the 
boundaries of $(D_{R}, g_{D})$ and $(E^{R}, g_{E})$ gives the estimate
\begin{equation} \label{e5.17}
|\Phi^{\widetilde g}(\widetilde g)| \leq  Ce^{-\sqrt{n}R}, 
\end{equation}
inside $U$, where $C$ is independent of $R$.

 This gives the construction of the approximate solutions 
$\widetilde g$. In fact the construction gives a smooth moduli space 
$\widetilde{\mathcal M}_{\varepsilon}$ of approximate solutions on $M$, 
naturally diffeomorphic to the moduli space ${\mathcal M}_{\varepsilon}$ on 
$\Sigma$. One has a natural boundary map 
\begin{equation} \label{e5.18}
\widetilde \Pi: \widetilde{\mathcal M}_{\varepsilon} \rightarrow{\mathcal C} , 
\end{equation}
which is the constant map to the conformal class $[\gamma_{0}(\beta)]$; note 
that dim $\widetilde{\mathcal M} = m = 6g-3$.

\medskip

 From (5.17), for $\varepsilon$ sufficiently small, or equivalently $R$ sufficiently 
large, one would expect that it should be possible to perturb the 
metrics $\widetilde g$ to exact Einstein metrics $g$ on $M$, i.e. perturb 
$\widetilde g$ to $g$ satisfying $\Phi^{\widetilde g}(g) = 0$. Following 
the proof of Theorem 2.1, it is not difficult to show that 
$\Phi^{\widetilde g}$ is a submersion at $\widetilde g$, so that the 
image of $\Phi^{\widetilde g}$ contains an open ball $B(\mu) \subset {\mathcal C}$ 
about $\Phi^{\widetilde g}(\widetilde g)$. However, one needs to obtain a lower 
bound on the radius $\mu$ of such a ball, independent of $R$, or at 
least prove that 
\begin{equation} \label{e5.19}
\mu  >> e^{-\sqrt{n}R}, 
\end{equation}
for $R$ large. 

 Now it is clear that the linearized operator $L$ in (2.8) cannot be 
uniformly invertible at $\widetilde g \in\widetilde{\mathcal M}$. In fact, 
there is an approximate kernel $\widetilde K$ of $L$ acting on forms in 
${\mathbb S}_{2}^{m,\alpha}$, induced by forms $\kappa$ tangent to the 
moduli space ${\mathcal M}$ of hyperbolic metrics on $\Sigma$ and extended 
to $M$ in a natural way, so that, 
if $||\widetilde \kappa||_{L^{\infty}} = 1$, then 
\begin{equation} \label{e5.20}
||L(\widetilde \kappa)||_{L^{\infty}} \rightarrow  0, 
\end{equation}
as the glueing radius $R \rightarrow  \infty$. Note that forms $\kappa$ 
tangent to ${\mathcal M}$, (equivalent to holomorphic quadratic 
differentials on $\Sigma$), decay to 0 on the end $E$ outside the closed 
geodesic $\sigma$. 

  Suppose that the Einstein manifold $(N, g_{N})$ in (2.10) has the 
property that 
\begin{equation} \label{e5.21}
K_{N} \leq  0, 
\end{equation}
i.e. $g_{N}$ has non-positive curvature. Then it is not particularly 
difficult to show, although we will not give the details here, that $L$ 
is uniformly invertible on the orthogonal complement 
$\widetilde K^{\perp}$ of $\widetilde K$ in ${\mathbb S}_{2}^{m,\alpha}$ with 
respect to the $L^{2}$ metric. One has then dim$\widetilde K = m = 6g-3$. 

 There are now two methods to try to establish (5.19) and thus obtain a 
Poincar\'e-Einstein metric $g$ near $\widetilde g$. First, one can try to 
arrange that 
\begin{equation} \label{e5.22}
\Phi^{\widetilde g}(\widetilde g) \in  \widetilde K^{\perp}, 
\end{equation}
possibly by modifying or perturbing $\widetilde g$, and iteratively solve 
$\Phi^{\widetilde g}(g) = 0$ within the space ${\mathbb S}_{2}^{m,\alpha}$. 
This would lead to the existence of a Poincar\'e-Einstein metric $g$ with 
the same boundary metric as $\widetilde g$. However this is not possible: 

\begin{proposition} \label{p 5.3.}
  There is no Einstein metric on $M$ with conformal infinity given by 
$S^{1}(\beta)\times (N, g_{N})$. 
\end{proposition}

\noindent
{\bf Proof:}
 This is an immediate consequence of Corollary 2.5.
{\endproof}

 Thus, to obtain an Einstein metric $g$ on $M$ close to $\widetilde g$ 
requires changing the boundary metric of $\widetilde g$, i.e. working 
outside the space ${\mathbb S}_{2}^{m,\alpha}$. This means one must use 
the dependence of $\Phi^{\widetilde g}$ on the boundary metrics to try 
to kill the cokernel of $L$, as in the proof of Theorem 2.1. In turn, this 
requires proving that the pairing (2.9) is non-degenerate and bounded below. 
More precisely, for any $\widetilde \kappa\in\widetilde K$, there must 
exist a variation $\dot \gamma$ of the boundary metric $\gamma_{0}$, with 
$||\widetilde \kappa||_{L^{\infty}} \leq 1$ and  
$||\dot \gamma||_{L^{\infty}} \leq 1$, such that 
\begin{equation} \label{e5.23}
|\int_{M}\langle D\Phi^{\widetilde g}(\dot \gamma), 
\widetilde \kappa \rangle dV| \geq  \mu' ,   
\end{equation}
 where $\mu'$ also satisfies (5.19). 

 We have not been able to verify (5.23), and expect it is not true. 
First, the choice of $\widetilde K$ is not canonical, and to establish 
(5.23), one needs a precise definition or choice. If $\widetilde K$ is defined 
as forms tangent to the moduli space ${\mathcal M}_{\varepsilon}$ on $\Sigma$, 
naturally extended to forms on $M$ and which have compact support, then 
$\kappa  = 0$ near infinity. However, by the definition given in \S 2, 
the support of $D\Phi^{\widetilde g}(\dot \gamma)$ is near 
$\partial M$. Thus, this definition implies that the integrals in 
(5.23) are all 0. 

 However one defines $\widetilde K$ exactly, (for instance as the space 
of eigenforms with small eigenvalues of $L$ acting on ${\mathbb 
S}_{2}^{m,\alpha}$), the support of any $\widetilde \kappa \in  
\widetilde K$ with $L^{\infty}$ norm 1 will be almost completely 
contained in the region $D$ from (5.14). Thus, such forms must decay 
quickly at infinity, and we expect the decay is too fast to give the 
lower bound (5.23). However, we have not been able to verify this in 
detail. Instead, we give other, heuristic, arguments which suggest that 
(5.23) must fail.

\medskip

 First, when (5.21) holds, a straightforward computation shows that the 
extremal black hole $(N, g_{E})$ also has non-positive curvature, cf. [10] 
for instance. Thus, (5.1) holds and Conjecture 5.1 would imply that any metric 
$\gamma$ on the boundary $S^{1}\times N$ sufficiently close to 
$\gamma_{0}(\beta)$ is the boundary metric of a Poincar\'e-Einstein 
$g_{\gamma}$ with a cusp end on the manifold ${\mathbb R}\times S^{1}\times N$. 
The metric $g_{\gamma}$ is asymptotic to the extreme metric $g_{E}$, or 
equivalently to the metric $g_{\infty}$ in (5.8) down the cusp end. This is 
because $g_{\infty}$ is rigid, in that it admits no bounded infinitesimal 
Einstein deformations. Consequently, exactly the same construction of the 
approximate Einstein metrics as above may be carried out with $g_{\gamma}$ 
in place of the extremal black hole metric $g_{E}$. 

  This now gives an infinite dimensional moduli space $\widetilde{\mathcal E}$ of 
approximate Einstein metrics and a corresponding boundary map 
\begin{equation} \label{e5.24}
\widetilde \Pi: \widetilde{\mathcal E} \rightarrow{\mathcal C} . 
\end{equation}
The map $\widetilde \Pi$ is surjective onto a neighborhood of 
$\gamma_{0}(\beta )$ and has fibers $\widetilde \Pi^{-1}(\gamma) = 
\widetilde{\mathcal M}_{\gamma}$ diffeomorphic to $\widetilde{\mathcal M}$ as 
before. Thus, $\widetilde \Pi$ is Fredholm of index $m = 6g-3$. Clearly 
$\widetilde \Pi$ is a submersion, so that the dimension of the kernel of 
$D_{\widetilde g}\widetilde \Pi$ is $m$, at every 
$\widetilde g\in\widetilde{\mathcal E}$. 

 Now if (5.23) holds, (for some choice of $\widetilde K$), one expects 
it should hold equally well on $\widetilde{\mathcal M}_{\gamma}$ in place 
of $\widetilde{\mathcal M}$. This implies the existence of a space of 
Poincar\'e-Einstein metrics ${\mathcal E}$ and smooth boundary map 
\begin{equation} \label{e5.25}
\Pi: {\mathcal E}  \rightarrow{\mathcal C} . 
\end{equation}
We abuse notation here slightly and assume that ${\mathcal E}$ consists 
only of metrics close to the approximate Einstein metrics 
$\widetilde g\in\widetilde{\mathcal E}$; thus ${\mathcal E}$ is a connected open 
subset of the full moduli space. 

  Both maps $\widetilde \Pi$ and $\Pi$ extend as continuous maps to the completions 
$\overline{\widetilde{\mathcal E}}$ and $\overline{\mathcal E}$ of the moduli spaces 
$\widetilde{\mathcal E}$ and ${\mathcal E}$ respectively. The boundary 
$\partial \widetilde {\mathcal E}$ consists of the (original) cusp metrics $g_{\gamma}$ 
constructed as perturbations of the extreme metric $g_{E}$; in particular, 
the metrics in $\partial \widetilde {\mathcal E}$ are all Einstein. For metrics 
$\widetilde g$ near $\partial \widetilde {\mathcal E}$, the estimate (5.17) holds, 
with $R$ very large. Hence, one has 
$$\partial \widetilde {\mathcal E} = \partial {\mathcal E},$$ 
and $\Pi$ maps $\partial {\mathcal E}$ onto an open set in ${\mathcal C}$. 
This of course contradicts Conjecture 4.2. In fact the maps 
$D_{\widetilde g}\widetilde \Pi$ and $D_{g}\Pi$ are both Fredholm, and are 
arbitrarily close for $R$ sufficiently large. The Fredholm index is constant 
under small perturbations and since index $D_{\widetilde g}\widetilde \Pi = 
6g-3$, while index $D_{\widetilde g}\widetilde \Pi = 0$, one has a contradiction. 

\medskip

 Although the arguments above require some further justifications to be 
made completely rigorous, they strongly suggest that (5.23) does not 
hold and that metrics in $\widetilde{\mathcal M}$ or $\widetilde{\mathcal E}$ 
cannot be perturbed to Poincar\'e-Einstein metrics on $M$. Further independent 
evidence for this is also given by Corollary 2.5, which implies that if the 
space ${\mathcal E}$, close to $\widetilde{\mathcal E}$ exists, then Im $\Pi$ 
must miss the infinite-dimensional space of $S^{1}$-invariant metrics 
on $\partial M$. We see no good reason why this should be the case. 

\medskip

 This leads again then to the basic question raised at the end of \S 4 
and discussed above; are there situations where Poincar\'e-Einstein cusps 
form as limits of metrics in a component ${\mathcal E}_{0}$ of ${\mathcal E}$? 
Resolution of this question would represent major progress on understanding the 
global existence question for the Dirichlet problem. 

\bibliographystyle{plain}

\begin{center}
March, 2005
\end{center}

\noindent
\address{Department of Mathematics\\
S.U.N.Y. at Stony Brook\\
Stony Brook, N.Y. 11794-3561}\\
\email{anderson@math.sunysb.edu}

\end{document}